\documentclass[11pt]{article}

\usepackage{fancybox,color,amsmath,amsfonts,amssymb,mathrsfs, graphicx,cases,url,bm}
\usepackage{geometry}
\usepackage{hyperref}

\definecolor{lred}{rgb}{1,0.8,0.8}
\definecolor{lblue}{rgb}{0.8,0.8,1}
\definecolor{dred}{rgb}{0.6,0,0}
\definecolor{dblue}{rgb}{0,0,0.5}
\definecolor{red}{rgb}{0.9,0,0}
\definecolor{blue}{rgb}{0,0,0.9}

\renewcommand{\Re}{{\rm I}\!  {\rm R}}

\newtheorem{proposition}{Proposition}[section]

\newtheorem{lemma}{Lemma}[section]

\newtheorem{definition}{Definition}[section]
\newtheorem{remark}{Remark}[section]

\newtheorem{claim}{Claim} 
\begin{document}

\title{Variational Analysis of the Ky Fan $k$-norm}

\author{Chao Ding\footnote{Institute of Applied Mathematics,  Chinese Academy of Sciences, Beijing,  P.R. China. Email: dingchao@amss.ac.cn. This work is supported
		in part by the National Natural Science Foundation
		of China (Grant No. 11301515).}}
\date{October 19, 2015}

%\date{Received: date / Accepted: date}
% The correct dates will be entered by the editor

\maketitle

\begin{abstract}
In this paper, we will study some variational properties of the Ky Fan $k$-norm $\theta=\|\cdot\|_{(k)}$ of matrices, which are closed related to a class of basic nonlinear optimization problems involving the Ky Fan $k$-norm. In particular, for the basic nonlinear optimization problems, we will introduce the concept of nondegeneracy, strict complementarity and the critical cones associated with  the generalized equations. Finally, we present the explicit formulas of the conjugate function of the parabolic second order directional derivative of $\theta$, which will be referred to as the {\it sigma term} of the second order optimality conditions.  The results obtain in this paper  provide the necessary theoretical foundations {for} future work on sensitivity and stability analysis of the nonlinear optimization problems involving the Ky Fan $k$-norm.
\end{abstract}

\vskip 10 true pt



%%%%%%%%%%%%%%%%%%%%%%%%%%%%%%%%%%%%%%%%%%%%%%%%%%%%%%%%%

\section{Introduction}\label{sect:Intro}

Let $\Re^{m\times n}$ be the vector space of all $m\times n$ real matrices equipped with the inner product $\langle Y,Z\rangle:={\rm Tr}(Y^TZ)$ for $Y$ and $Z$ in $ \Re^{m\times n}$, where $``{\rm Tr}"$ denotes the trace, i.e., the sum of the diagonal entries,  of a squared matrix. For simplicity, we always assume that $m\le n$.  For any given positive integer $1\le k\le m$, denote $\theta:=\|\cdot\|_{(k)}$ the matrix Ky Fan $k$-norm, i.e., the sum of $k$ largest singular values of matrices. In particular, $\|\cdot\|_{(1)}$ coincides with the spectral norm $\|\cdot\|_2$ of the matrices, i.e., the largest singular value of matrices; $\|\cdot\|_{(m)}$ is the nuclear norm $\|\cdot\|_*$ of matrices, i.e., the sum of singular values of matrices. It is well-known that $\vartheta(Z):=\|Z\|^*_{(k)}=\max\{\|Z\|_2,\|Z\|_*/k\}$ for $Z\in\Re^{m\times n}$ is the dual norm of $\|\cdot\|_{(k)}$ (cf. \cite[Exercise IV.1.18]{Bhatia97}).  Since $\theta$ is a matrix norm (convex, closed, positively homogeneous and $\theta(0)=0$), we obtain from \cite[Theorem 13.5 \& 13.2]{Rockafellar70} that the conjugate function $\theta^*=\delta_{\partial\,\theta(0)}$ is just the indicator function of the subdifferential $\partial\,\theta(0)$ of $\theta$ at $0$.  Moreover, it can be verified directly from the definition of dual norm that $\partial\,\theta(0)$ coincides with the unit ball under the dual norm $\vartheta$, i.e., $\partial\,\theta(0)={\mathbb B}_{(k)^*}:=\{S\in\Re^{m\times n}\mid \vartheta(S)\leq 1 \}$.

Consider the following nonlinear optimization problem involving the Ky Fan $k$-norm $\theta=\|\cdot\|_{(k)}$ 
\begin{equation}\label{eq:P-g}
\min\left\{ f(x)+\theta(g(x))\mid x\in{\cal X}\right\},
\end{equation}
where ${\cal X}$ is a finite dimensional real vector space equipped with a scalar product $\langle\cdot,\cdot\rangle$,  $f:{\cal X}\to\Re$ is a continuously differentiable real value function,  and $g:{\cal X}\to\Re^{m\times n}$ is a continuously differentiable function. Since $\theta$ is convex and finite everywhere, it is well-known \cite[Example 10.8]{RWets98} that for a locally optimal solution $\bar{x}\in{\cal X}$ of \eqref{eq:P-g}, there always exists a Lagrange multiplier $\overline{S}\in\Re^{m\times n}$, together with $\bar{x}$ satisfying the following first order optimality condition, namely the Karush-Kuhn-Tucker (KKT) condition:
\begin{equation}\label{eq:KKT-g}
\nabla f(\bar{x})+g'(\bar{x})^*\overline{S}=0 \quad {\rm and} \quad \overline{S}\in\partial\,\theta(\overline{X}),
\end{equation}
where $\overline{X}:=g(\bar{x})$, $\nabla f(\bar{x})\in{\cal X}$ is the gradient of $f$ at $\bar{x}$, $g'(\bar{x})^*:\Re^{m\times n}\to{\cal X}$ is the adjoint of the derivative mapping $g'(\bar{x})$. Note that if ${\cal X}=\Re^{m\times n}$ and $g(x):=x$ is the identity mapping, then the KKT condition \eqref{eq:KKT-g} becomes the following generalized equation:
\[
0\in \nabla f(x)+\partial\,\theta(x).
\]
Note also that if the function $\theta$ is replaced by the indicator function $\delta_{K}$ of a set $K$ in a finite dimensional real vector space, then the nonlinear optimization problem \eqref{eq:P-g} becomes 
\begin{equation}\label{eq:simple-P}
	\begin{array}{cl} \min &  f(x)
\\ [3pt]
{\rm s.t.} & g(x)\in K.
\end{array}
\end{equation} 
During the last three decades, considerable progress has been made in the variational analysis related to the problem \eqref{eq:simple-P} \cite{RWets98,FPang03,BShapiro00,KKummer02,Mordukhovich06a}. In particular, for the general non-polyhedral set $K$ (e.g., the second-order cone and the positively semidefinite (SDP) matrices cone), by employing the well studied properties of the variational inequality $S\in{\cal N}_{K}(x)$, some important  properties of \eqref{eq:simple-P}, such as the constraint nondegeneracy, second order optimality conditions, strong regularity, full stability and calmness, are studied recently by various researchers \cite{BShapiro00,Sun06,MRSarabi13}. In order to extend those results to the optimization problems  involving the Ky Fan $k$-norm, we need first study the variational properties of \eqref{eq:P-g}, especially the properties of the generalized equation $S\in\theta(X)$ and its equivalent dual problem $X\in\theta^*(S)$. Although the optimization problem \eqref{eq:P-g} seems extremely simple, many fundamental and important issues such as the concept of nondegeneracy, the characterizations of critical cones and the second order optimality conditions, are not studied yet in literature. The main purpose of this paper is to build up the necessary variational foundations for the future work on the nonlinear optimization problems involving the Ky Fan $k$-norm. 

Certainly, instead of the {\it basic} model \eqref{eq:P-g}, one can consider its various modifications, e.g., the nonlinear optimization problems involving the Ky Fan $k$-norm with equality and conic constraints. In particular, the following convex composite matrix optimization problems involving the Ky Fan $k$-norm frequently arise in various applications such as the matrix norm approximation, matrix completion, rank minimization, graph theory, machine learning, etc \cite{GT94,Toh97,TT98,RFParrilo10,CanRec08,CanTao09,CLMW09,CSPW09,WMGR09,CFP03,BDXiao04,GSun10,Lovasz79,Dobrynin04,KLVempala97,MSPan12}:
\begin{equation}\label{eq:def-MOP}
\begin{array}{cl} \min &  \frac{1}{2}\langle (X,Y),{\cal Q}(X,Y)\rangle+\langle C,(X,Y)\rangle +  \theta(X) 
\\ [3pt]
{\rm s.t.} & {\cal A}(X,Y) = b, \quad Y\in K,
\end{array}
\end{equation}
where  ${\cal Y}$ is a  finite dimensional real vector space, ${\cal Q}:\Re^{m\times n}\times {\cal Y}\to \Re^{m\times n}\times{\cal Y}$ is a positively semidefinite self-adjoint linear operator, ${\cal A}:\Re^{m\times n}\times{\cal Y}\to\Re^p$ is a linear operator, $C\in\Re^{m\times n}\times{\cal Y}$ and $b\in\Re^p$ are given data, and $K\in{\cal Y}$ is a closed convex cone (e.g., the positive orthant, second-order cone of vectors, positive semidefinite matrices cone). As the initial step, in this paper, we will mainly focus on the fundamental model \eqref{eq:P-g}, since the obtained variational results will provide the necessary theoretical foundations for the study of more complicate model, e.g., \eqref{eq:def-MOP}. More precisely, we will study the concepts of nondegeneracy and strict complementary to locally optimal solutions of \eqref{eq:P-g}. Also, we will define and provide the complete characterizations of the critical cones associated with  the generalized equation $S\in\theta(X)$ and its dual problem $X\in\theta^*(S)$.  Another important variational property studied in this paper is the conjugate function of the parabolic second order directional derivative of the Ky Fan $k$-norm $\theta$, which equals to the support function of the second order tangent set of the epigraph of $\theta$.  This conjugate function is closely related to the second order optimality conditions of the problem \eqref{eq:P-g}. Note that the epigraph of $\theta$ is not polyhedral. In general, the conjugate function of the parabolic second order directional derivative of the Ky Fan $k$-norm $\theta$ will not vanish in the corresponding second order optimality conditions, and will be referred to as the {\it sigma term}, provides the second order information of $\theta$. In this paper, we provide the explicit expression of this sigma term. Consequently, it becomes possible to establish the  second order optimality conditions of the problem \eqref{eq:P-g} and study many corresponding sensitivity properties, e.g., the second order optimality conditions and the characterization of strong regularity of the KKT solutions.

The remaining parts of this paper are organized as follows. In Section \ref{sect:preliminary}, we introduce some preliminary results on the differential properties of eigenvalue values and vectors of symmetric matrices and singular values and vectors of matrices. In Section \ref{sect:solutions of GEs}, we study the properties of the solution of the GE $S\in\partial\,\theta(X)$, which arises from the KKT condition \eqref{eq:KKT-g} and its equivalent dual form $X\in\partial\,\theta^*(S)$. We introduce the nondegeneracy and strict complementarity of \eqref{eq:P-g} in Section \ref{sect:nondegeneracy}. In Section \ref{sect:critical cones}, we introduce and study the critical cones associated with the GE $S\in\partial\,\theta(X)$ and $X\in\partial\,\theta^*(S)$. The second order properties of the Ky Fan $k$-norm $\theta$ are studied in Section \ref{sect:second order tangent sets}. We conclude our paper in the final section.
 
Below are some common notations to be used:
\begin{itemize}
	\item For any $Z\in\Re^{m\times n}$, we denote by $Z_{ij}$ the $(i,j)$-th entry of $Z$.
	\item For any    $Z\in\Re^{m\times n}$,  we use  $z_{j}$ to represent the $j$th column of $Z$,
	$j=1,\ldots,n$.  Let ${\cal J}\subseteq  \{1,\ldots, n\}$ be an index set. We use
	$ Z_{{\cal J}}$ to  denote the sub-matrix of $Z$ obtained by  removing all the columns of $Z$ not in
	${\cal J}$.
	So for each $j$, we have   $ Z_{\{j\}}=z_j$.
	
	\item   Let ${\cal I}\subseteq \{1,\ldots, m\}$ and  ${\cal J}\subseteq \{1,\ldots, n\}$ be two index
	sets. For  any    $Z\in\Re^{m\times n}$, we use $Z_{{\cal I}{\cal J}}$ to
	denote the $|{\cal I}|\times|{\cal J}|$ sub-matrix of $Z$ obtained  by removing all the rows of $Z$
	not in ${\cal I}$ and all the columns of $Z$ not in  ${\cal J}$.
	\item We use $``\circ"$ to denote the Hardamard product between matrices, i.e., for any two
	matrices $X$ and $Y$ in $\Re^{m\times n}$ the $(i,j)$-th entry of $  Z:= X\circ Y \in
	\Re^{m\times n}$ is
	$Z_{ij}=X_{ij} Y_{ij}$.
\end{itemize}

\section{Preliminaries}\label{sect:preliminary}

In this section, we list some useful preliminary results on the eigenvalues of   symmetric matrices and the singular values of matrices, which are useful for our subsequent analysis.

Let ${\cal  S}^{n}$ be the space of all real $n\times n$ symmetric matrices and ${\cal  O}^{n}$ be the
set of all $n\times n$ orthogonal matrices.  Let  $X \in {\cal  S}^n$ be given.  We use
$\lambda_{1}(X)\ge \lambda_2(X) \ge \ldots \ge \lambda_{n}(X)$ to denote the real eigenvalues of $X$
(counting multiplicity) being arranged in non-increasing order.   Denote
$\lambda(X):=(\lambda_1(X),\lambda_2(X), \ldots, \lambda_n(X))^{T}\in\Re^{n}$ and $\Lambda(X):=
{\rm diag}(\lambda(X))$, where for any $x\in \Re^n$, ${\rm diag}(x)$ denotes the diagonal matrix
whose $i$-th diagonal entry is $x_i$, $i=1, \ldots, n$.
Let  $P\in{\cal  O}^{n}$ be such that
\begin{equation}\label{eq:eig-decomp}
X=P\Lambda(X)P^{T} \, .
\end{equation}
 We denote the set of such matrices $P$ in the eigenvalue decomposition (\ref{eq:eig-decomp})  by
 ${\cal  O}^{n}(X)$. Let $\omega_{1}(X)>\omega_{2}(X)>\ldots>\omega_{r}(X)$ be the distinct eigenvalues of $X$. Define the index sets
 \begin{equation}\label{eq:ak-symmetric}
a_{k} :=\{i\,|\, \lambda_{i}(X)={\omega}_{k}(X),\ 1\leq i\leq n\},\quad k=1,\ldots,r.
 \end{equation} For each $  i\in \{1,\ldots,n\}$, we define $l_{i}(X)$ to be the number of eigenvalues
 that are equal to $\lambda_i(X)$ but are ranked  before $i$ (including $i$) and $s_{i}(X)$ to be the
 number of eigenvalues that are equal to $\lambda_i(X)$ but are ranked  after $i$ (excluding $i$),
 respectively, i.e.,
 we define $l_{i}(X)$ and  $s_{i}(X)$  such that
\begin{eqnarray}
&&\lambda_{1}(X)\geq\ldots\geq\lambda_{i-l_{i}(X)}(X)>\lambda_{i-l_{i}(X)+1}(X)=\ldots=\lambda_{i}(X)=\ldots=\lambda_{i+s_{i}(X)}(X)\nonumber\\ 
&&>\lambda_{i+s_{i}(X)+1}(X)\geq\ldots\geq\lambda_{n}(X).\label{eq:symmetric-l_i}
\end{eqnarray}
In later discussions, when the dependence of $l_{i}$ and $s_{i}$, $i=1, \ldots, n$,  on $X$ can be
seen clearly from the context, we often drop $X$ from these notations.

Next, we list some useful results about the symmetric matrices which are needed in subsequent
discussions. The inequality in the following lemma is known as Fan's inequality \cite{Fan49}.
 \begin{lemma}\label{lem:Fan}
Let   $Y$ and $Z$ be  two matrices in ${\cal  S}^{n}$. Then
\begin{equation}\label{eq:Fan}
\langle Y,Z\rangle\leq\lambda(Y)^{T}\lambda(Z)\, .
\end{equation}
where the equality holds if and only if $Y$ and $Z$ admit a simultaneous ordered eigenvalue
decomposition, i.e., there exists an orthogonal matrix $U\in{\cal  O}^{n}$ such that
\[
Y=U\Lambda(Y)U^{T}\quad {\rm and} \quad Z=U\Lambda(Z)U^{T}.
\]
\end{lemma}

\vskip 10 true pt

The following proposition on the directional differentiability of the eigenvalue function
$\lambda(\cdot)$ is well known. For example, see \cite[Theorem 7]{Lancaster64} and \cite[Proposition
1.4]{Torki01}.

 \begin{proposition}\label{prop:eigenvalue_diff}
Let $X\in {\cal  S}^{n}$   have the eigenvalue decomposition (\ref{eq:eig-decomp}).
Then, for any ${\cal  S}^n\ni H\to 0$, we have
\begin{equation}\label{eq:dirctdiff_eigen_genera}
\lambda_{i}(X+H)-\lambda_{i}(X)-\lambda_{l_{i}}({P}_{a_{k}}^{T}H{P}_{a_{k}})=O(\|H\|^{2}),
\quad i\in \alpha_{k}, \ k=1,\ldots,r,
\end{equation}
where for each $i\in\{1,\ldots,n\}$, $l_{i}$ is defined in (\ref{eq:symmetric-l_i}).
Hence,  for any given direction $H\in {\cal  S}^{n}$, the eigenvalue function $\lambda_{i}(\cdot)$ is
directionally differentiable at $X$  with
$\lambda'_{i}(X;H)=\lambda_{l_{i}}({P}_{a_{k}}^{T}H{P}_{a_{k}})$, $i\in a_{k}$, $k=1,
\ldots, r$.
\end{proposition}

Let $k\in\{1,\ldots,r\}$ be fixed. For the symmetric matrix $P_{a_{k}}^{T}HP_{a_{k}}\in{\cal
S}^{|a_{k}|}$, consider the eigenvalue decomposition
\begin{equation}\label{eq:decop_PHP}
P_{a_{k}}^{T}HP_{a_{k}}=R\Lambda(P_{a_{k}}^{T}HP_{a_{k}})R^{T},
\end{equation}
where $R\in{\cal  O}^{|a_{k}|}$. Denote the distinct eigenvalues of $P_{a_{k}}^{T}HP_{a_{k}}$ by
$\tilde{\mu}_{1}>\tilde{\mu}_{2}>\ldots>\tilde{\mu}_{\tilde{r}}$. Define
\begin{equation}\label{eq:a-tilde-k-symmetric}
\tilde{a}_{j}:=\{i\,|\, \lambda_{i}(P_{a_{k}}^{T}HP_{a_{k}})=\tilde{\mu}_{j}, 1\leq i\leq
|a_{k}|\},\quad j=1,\ldots,\tilde{r}.
\end{equation} For each $i\in a_{k}$, let $\tilde{l}_{i}\in\{1,\ldots,|a_{k}|\}$ and
$\tilde{k}\in\{1,\ldots,\tilde{r}\}$ be such that
\begin{equation}\label{eq:tilde-li-k}
\tilde{l}_{i}:=l_{l_{i}}(P_{a_{k}}^{T}HP_{a_{k}}) \quad {\rm and} \quad
\tilde{l}_{i}\in\tilde{a}_{\tilde{k}},
\end{equation}  where $l_{i}$ is defined by (\ref{eq:symmetric-l_i}).

Let ${\cal X}$ and ${\cal X}'$ be two finite dimensional real Euclidean spaces. We say that a function $\Phi:{\cal X}\to{\cal X}'$ is (parabolic) second order directionally differentiable at $x\in{\cal X}$, if $\Phi$ is directionally differentiable at $x$ and for any $h,w\in{\cal X}$,
\[
\lim_{t\downarrow 0}\frac{\Phi(x+th+\frac{1}{2}t^2w)-\Phi(x)-t\Phi(x;h)}{\frac{1}{2}t^2}\quad \mbox{exists;}
\]
and the above limit is said to be the (parabolic) second order directional derivative of $\Phi$ at $x$ along the directions $h$ and $w$, denoted by $\Phi''(x;h,w)$. The following proposition \cite[Proposition 2.2]{Torki01}, provides the explicit formula of the (parabolic)
second order directional derivative of the eigenvalue function.

\begin{proposition}\label{prop:second-directional-diff-eigenvalue}
Let $X\in{\cal  S}^{n}$ have the eigenvalue decomposition (\ref{eq:eig-decomp}). Then, for any given
$H,W\in{\cal  S}^{n}$, we have for each $k\in\{1,\ldots,r\}$
\begin{equation}\label{eq:second-directional-diff-eigenvalue}
\lambda''_{i}(X;H,W)=\lambda_{\tilde{l}_{i}}\left(R_{\tilde{a}_{\tilde{k}}}^{T}P_{a_{k}}^{T}\left[W-2H(X-\lambda_{i}I_{n})^{\dag}H\right]P_{a_{k}}R_{\tilde{a}_{\tilde{k}}}\right),
\quad i\in a_{k},
\end{equation}
where $Z^{\dag}\in \Re^{n\times n}$ is the Moore-Penrose pseudoinverse of the square matrix $Z\in\Re^{n\times n}$.
\end{proposition}

Let  $X\in\Re^{m\times n}$ be given. Without loss of generality, assume that $m\leq n$. We use
$\sigma_{1}(X)\geq\sigma_{2}(X)\geq\ldots\geq\sigma_{m}(X)$ to
denote the singular values of $X$ (counting multiplicity) being
arranged in non-increasing order. Define
$\sigma(X):=(\sigma_{1}(X),\sigma_{2}(X),\ldots,\sigma_{m}(X))^{T}$
and $\Sigma(X):={\rm diag}(\sigma(X))$. Let $X\in\Re^{m\times n}$
admit the following singular value decomposition (SVD):
\begin{equation}\label{eq:SVD}
X={U}\left[ \Sigma(X)\ \  0 \right]{V}^{T}={U}\left[ \Sigma(X)\ \  0\right]\left[ {V}_{1} \ \  {V}_{2}
\right]^{T}={U}\Sigma(X) {V}_{1}^{T},
\end{equation} where
 ${U}\in{\cal  O}^{m}$ and ${V}=\left[ {V}_{1} \ \  {V}_{2} \right] \in{\cal  O}^{n}$ with
  ${V}_{1}\in\Re^{n\times m}$ and
  ${V}_{2}\in\Re^{n\times (n-m)}$.
   The set of such matrices pair $(U,V)$ in the SVD (\ref{eq:SVD}) is denoted by ${\cal  O}^{m,n}(X)$,
   i.e.,
\[
{\cal  O}^{m,n}(X):=\left\{(U,V)\in{\cal  O}^{m}\times{\cal  O}^{n}\,|\, X=U\left[ \Sigma(X)\ \  0
\right]V^{T}\right\}.
\] Define the three index sets $a$, $b$ and $c$ by
\begin{equation}\label{eq:def-a-b-c}
a:=\{i \,|\,\sigma_{i}(X)>0,\ 1\le i\le m\}, \  b:=\{i\,|\,\sigma_{i}(X)=0, \ 1\le i\le m\} \  {\rm and} \
c:=\{m+1,\ldots,n\}.
\end{equation}
 Let ${\nu}_{1}(X)>{\nu}_{2}(X)>\ldots>{\nu}_{r}(X)>0$ be the distinct nonzero singular values of $X$. Without causing any ambiguity, we also use $a_k$ to denote the following index sets
\begin{equation}\label{eq:ak-nonsymmetric}
a_{k}:=\{i\,|\,\sigma_{i}(X)={\nu}_{k}(X), \ 1\le i\le m\}, \quad k=1,\ldots,r.
\end{equation} For the sake of convenience, let $a_{r+1}:=b$.
For each $  i\in \{1,\ldots,m\}$, we also define $l_{i}(X)$ to be the number of singular values that are
equal to $\sigma_ i(X)$ but are ranked  before $i$ (including $i$) and $s_{i}(X)$ to be the number of
singular values that are equal to $\sigma_i(X)$ but are ranked  after $i$ (excluding $i$), respectively,
i.e.,
 we define $l_{i}(X)$ and  $s_{i}(X)$  such that
\begin{eqnarray}
&&\sigma_{1}(X)\geq\ldots\geq\sigma_{i-l_{i}(X)}(X)>\sigma_{i-l_{i}(X)+1}(X)=\ldots=\sigma_{i}(X)=\ldots=\sigma_{i+s_{i}(X)}(X)\nonumber\\ 
&&>\sigma_{i+s_{i}(X)+1}(X)\geq\ldots\geq\sigma_{m}(X).\label{eq:onsymmetric-l_i}
\end{eqnarray}
In later discussions, when the dependence of $l_{i}$ and $s_{i}$, $i=1, \ldots, m$,  on $X$ can be
seen clearly from the context, we often drop $X$ from these notations. The inequality in the following lemma is known as von Neumann's trace inequality \cite{vonNeumann37}.
\begin{lemma}\label{lem:vonNeumann}
	Let $Y$ and $Z$ be two matrices in $\Re^{m\times n}$. Then
	\begin{equation}\label{eq:vonNeumann}
	\langle Y, Z\rangle \le \sigma(Y)^T \sigma(Z),
	\end{equation} where the equality holds if $Y$ and $Z$ admit a simultaneous ordered singular value decomposition, i.e., there exist orthogonal matrices $U\in{\cal O}^{m}$ and $V\in{\cal O}^{n}$ such that 
	\[
	Y=U[\Sigma(Y)\ \ 0]V^{T}\quad {\rm and} \quad Z=U[\Sigma(Z)\ \ 0]V^{T}.
	\]
\end{lemma}

%Let ${\cal
%B}(\cdot):\Re^{m\times n}\to{\cal  S}^{m+n}$ be the linear
%operator defined by
%\begin{equation}\label{eq:B}
%{\cal B}(Z): =\left[ \begin{array}{cc} 0 & Z
%\\
%Z^T & 0
%\end{array} \right] , \quad Z\in \Re^{m\times n}.
%\end{equation}  We use $I^{\uparrow}_p$ to denote the $p$ by $p$
%anti-diagonal matrix whose anti-diagonal entries are all ones and
%other entries are zeros. Denote
%\[
%{U}_{a}^{\uparrow}={U}_{a} I^{\uparrow}_{|a|} \quad {\rm and} \quad
%{V}_{a}^{\uparrow}={V}_{a}I^{\uparrow}_{|a|} .
%\]  Let
%\begin{equation}\label{eq:P-def}
%{P}:=\frac{1}{\sqrt{2}}\left[\begin{array}{cccccc}{U}_{a} &
%{U}_{b} & 0 & {U}_{b} & {U}_{a}^{\uparrow} \\ {V}_{a} & {V}_{b} &
%\sqrt{2}\,{V}_{2} & -{V}_{b} & -{V}_{a}^{\uparrow}
%\end{array}\right]\in{\cal  O}^{m+n}.
%\end{equation} It is well-known  \cite[Theorem 7.3.7]{HJohnson85} that
%\begin{equation}\label{eq:decom-B}
%{{P}}^{T}{\cal B}(X) {P} =\Lambda({\cal B}(X))=\left[\begin{array}{ccc} \Sigma(X)  & 0 & 0 \\0  & 0 & 0
%\\ 0 & 0 & -\Sigma(X)^{\uparrow}\end{array}\right],
%\end{equation}
%where $\Sigma(X)^{\uparrow}:=\Sigma(X)I^{\uparrow}_{m}$. 
For notational convenience, define two linear  operators
$S:\Re^{p\times p}\to {\cal  S}^{p}$ and $T:\Re^{p\times p}\to
\Re^{p\times p}$  by
\begin{equation}\label{eq:maps-ST}
S(Z):=\frac{1}{2}(Z+Z^{T}) \quad {\rm and} \quad
T(Z):=\frac{1}{2}(Z-Z^{T}) \quad \forall\, Z\in\Re^{p\times p}.
\end{equation} 
The following proposition on the directional derivatives of the singular value functions can be obtained directly from Proposition \ref{prop:eigenvalue_diff}. For
more details, see \cite[Section 5.1]{LSendov05} .

\begin{proposition}\label{prop:easy}
Let $X\in\Re^{m\times n}$ have the singular value decomposition (\ref{eq:SVD}). For any $\Re^{m\times
n}\ni H\to 0$, we have
\begin{equation}\label{eq:directi-diff}
\sigma_{i}(X+H)-\sigma_{i}(X)-\sigma'_{i}(X;H)=O(\|H\|^{2}),  \quad i=1,\ldots,m,
\end{equation} with
\begin{equation}\label{eq:directi-diff-sigma}
\sigma'_{i}(X;H)=\left\{ \begin{array}{lcl}
\lambda_{l_{i}}\left( S({U}_{a_{k}}^{T}H {V}_{a_{k}})\right) & {\rm if} & i\in a_{k},\ k=1,\ldots,r,\\ 
\sigma_{l_{i}}\Big(\left[{U}_{b}^{T}H {V}_{b}\ \ {U}_{b}^{T}H {V}_{2} \right]\Big) & {\rm if} & i\in b,
\end{array}\right.
\end{equation}
where for each $i\in\{1,\ldots,m\}$, $l_{i}$ is defined in (\ref{eq:onsymmetric-l_i}).
\end{proposition}

Similarly, one can derive the following explicit formulas of the (parabolic) second order directional derivatives of the singular value functions from Proposition \ref{prop:second-directional-diff-eigenvalue}, directly. For more details, see \cite[Theorem 3.1]{ZZXiao13}.

\begin{proposition}\label{prop:second-directional-diff-singlevalue}
Let $X\in\Re^{m\times n}$ have the singular value decomposition (\ref{eq:SVD}). Suppose that the
direction $H,W\in \Re^{m\times n}$ are given. 
\begin{itemize}
\item[(i)] If $\sigma_{i}(X)>0$, then
\[
\sigma''_{i}(X;H,W)=\lambda_{\tilde{l}_{i}}\left( R_{\tilde{\alpha}_{\tilde{k}}}^{T} \left(
S(U^T_{a_k}WV_{a_k})-2 \Omega_{a_k}(X,H) \right)
R_{\tilde{\alpha}_{\tilde{k}}}\right),
\] 
where $k\in\{1,\ldots,r\}$ such that $i\in a_k$, $\Omega_{a_k}(X,H)\in{\cal S}^{m}$ is given by
\begin{eqnarray}
\Omega_{a_k}(X,H)&=&(S(U^THV_1)_{a_k})^T(\Sigma(X)-\nu_k(X)I_m)^{\dag}S(U^THV_1)_{a_k}\nonumber \\ [3pt]
&&+(T(U^THV_1)_{a_k})^T(-\Sigma(X)-\nu_k(X)I_m)^{\dag}T(U^THV_1)_{a_k} \nonumber \\ [3pt]
&& +\frac{1}{2\nu_k(X)}U^T_{a_k}HV_2V_2^TH^TU_{a_k},  \label{eq:def-Omega-dd}
\end{eqnarray}
%$\Omega_{k}(X,H)=\Gamma_1(X,S(U^THV_1)_{a_k})+\Gamma_2(X,T(U^THV_1)_{a_k})+\Gamma_3(X,U^T_{a_k}HV_2)$ the matrices $\Gamma_1(X,H)$,  $\Gamma_2(X,H)$ and $\Gamma_3(X,H)\in{\cal S}^{|a_k|}$ are given by
%\[
%\left\{
%\begin{array}{l}
%\Gamma_1(X,H) := (S(U^THV_1)_{a_k})^T(\Sigma(X)-\sigma_i(X)I_m)^{\dag}S(U^THV_1)_{a_k}, \\ [3pt]
%\Gamma_2(X,H):= (T(U^THV_1)_{a_k})^T(-\Sigma(X)-\sigma_i(X)I_m)^{\dag}T(U^THV_1)_{a_k}, \\ [3pt]
%\Gamma_3(X,H):= \frac{1}{2\sigma_i(X)}U^T_{a_k}HV_2(U^T_{a_k}HV_2)^T,
%\end{array}
%\right.
%\]
the matrix $R\in{\cal  O}^{|a_{k}|}$ satisfies $S(U_{a_{k}}^{T}HV_{a_{k}})=R\Lambda(S(U_{a_{k}}^{T}HV_{a_{k}}))R^{T}$,
and $\{\tilde{\alpha}_{j}\}_{j=1}^{\tilde{r}}$ and $\tilde{l}_{i}$, $\tilde{k}$  be defined by
(\ref{eq:a-tilde-k-symmetric}) and  (\ref{eq:tilde-li-k}) respectively for
$S(U_{a_{k}}^{T}HV_{a_{k}})$.
\item[(ii)] If $\sigma_{i}(X)=0$ and $\sigma_{l_{i}}([U_{b}^{T}HV_{b}\ \ U_{b}^{T}HV_{2}])>0$, then
\[
\sigma''_{i}(X;H,W)=\lambda_{\tilde{l}_{i}}(S(E_{\tilde{a}_{\tilde{k}}}^{T}[U_{b}^{T}ZV_{b}\ \
U_{b}^{T}ZV_{2}]F_{\tilde{a}_{\tilde{k}}})),
\] where $Z=W-2HX^{\dag}H\in\Re^{m\times n}$, $X^{\dag}\in\Re^{n\times m}$ is the Moore-Penrose pseudoinverse of  $X\in\Re^{m\times n}$, $E\in{\cal  O}^{|b|}$, $F=[F_{1}\ \ F_{2}]\in{\cal  O}^{|b|+(n-m)}$  satisfy
\[
[U_{b}^{T}HV_{b}\ \ U_{b}^{T}HV_{2}]=E[\Sigma([U_{b}^{T}HV_{b}\ \ U_{b}^{T}HV_{2}])\ \ 0]F^{T},
\] $\tilde{l}_{i}\in\{1,\ldots,|\tilde{a}_{\tilde{k}}| \}$ and $\tilde{k}\in\{1,\ldots,\tilde{r} \}$ such that $\tilde{l}_{i}=l_{l_i}(S(E_{\tilde{a}_{\tilde{k}}}^{T}[U_{b}^{T}ZV_{b}\ \
U_{b}^{T}ZV_{2}]F_{\tilde{a}_{\tilde{k}}}))$ and $\tilde{l}_{i}\in \tilde{a}_{\tilde{k}}$, $\tilde{a}_{j}$, $j=1,\ldots,\tilde{r}$ are the index sets of $[U_{b}^{T}HV_{b}\ \ U_{b}^{T}HV_{2}]$ defined by 
\[
\tilde{a}_{j}:=\{ i\,|\, \sigma_{i}([U_{b}^{T}HV_{b}\ \ U_{b}^{T}HV_{2}])=\tilde{\nu}_{j}, \ 1\leq i \leq |b|\},
\]
and $\tilde{\nu}_{1}>\tilde{\nu}_{2}>\ldots>\tilde{\nu}_{\tilde{r}}$ are the nonzero distinct singular values of $[U_{b}^{T}HV_{b}\ \ U_{b}^{T}HV_{2}]$.

\item[(iii)] If $\sigma_{i}(X)=0$ and $\sigma_{l_{i}}([U_{b}^{T}HV_{b}\ \ U_{b}^{T}HV_{2}])=0$, then
\[
\sigma''_{i}(X;H,W)=\sigma_{\tilde{l}_{i}}\Big( E_{\tilde{b}}^{T}[U_{b}^{T}ZV_{b}\
\ U_{b}^{T}ZV_{2}][F_{\tilde{b}}\ \ F_{2}] \Big),
\] where $Z=W-2HX^{\dag}H\in\Re^{m\times n}$, $\tilde{b}:=\{ i\,|\, \sigma_{i}([U_{b}^{T}HV_{b}\ \ U_{b}^{T}HV_{2}])=0, \ 1\leq i \leq |b|\}$ and $\tilde{l}_{i}=l_{l_{i}}\Big( E_{\tilde{b}}^{T}[U_{b}^{T}ZV_{b}\
\ U_{b}^{T}ZV_{2}][F_{\tilde{b}}\ \ F_{2}] \Big)$ is defined by
(\ref{eq:onsymmetric-l_i}) with respect to $E_{\tilde{b}}^{T}[U_{b}^{T}ZV_{b}\
\ U_{b}^{T}ZV_{2}][F_{\tilde{b}}\ \ F_{2}]$.

\end{itemize}
\end{proposition}

\section{The generalized equations}\label{sect:solutions of GEs}
In this section, we first study some properties of the following simple generalized equation (GE) 
\begin{equation}\label{eq:GE}
0\in -S+\partial\,\theta(X),
\end{equation} 
which is equivalent to the following dual form
\begin{equation}\label{eq:GE-dual}
0\in -X+\partial\,\theta^*(S).
\end{equation}  
Since $\theta^*=\delta_{{\mathbb B}_{(k)^*}}$, it follows from \cite[Theorem 23.5]{Rockafellar70} that \eqref{eq:GE} and \eqref{eq:GE-dual} are also equivalent to the following complementarity problem
\begin{equation}\label{eq:GE=CP}
(\overline{X},\theta(\overline{X}))\in{\cal K},\quad (\overline{S},-1)\in {\cal K}^{\circ} \quad {\rm and} \quad \left\langle(\overline{X},\theta(\overline{X})),(\overline{S},-1)\right\rangle=0,
\end{equation}
where ${\cal K}$ is the epigraph of  $\theta=\|\cdot\|_{(k)}$, i.e., 
\begin{equation}\label{eq:K-epsilon}
{\cal K}={\rm epi}\,\theta=\left\{(X,t)\in\Re^{m\times n}\times\Re \mid t\ge\|X\|_{(k)}\right\} 
\end{equation} 
and ${\cal K}^{\circ}$ is the polar cone of ${\cal K}$ given by 
\[
{\cal K}^{\circ}=\bigcup_{\rho\ge 0}\rho(\partial \theta(0),-1)=-{\rm epi}\,\vartheta \quad \mbox{with $\vartheta=\|\cdot\|_{(k)}^*$.}
\] 

On the other hand, it is well-known \cite{Moreau65} (see also \cite[Theorem 31.5]{Rockafellar70}) that $(\overline{X},\overline{S})$ is a solution of the GE \eqref{eq:GE} (or \eqref{eq:GE-dual}) if and only if
\[
\overline{X}-{\rm Pr}_{\theta}(\overline{X}+\overline{S})=0 \quad \Longleftrightarrow \quad \overline{S}-{\rm Pr}_{\theta^*}(\overline{X}+\overline{S})=0,
\]
where ${\rm Pr}_{\theta}:\Re^{m\times n}\to\Re^{m\times n}$ is the Moreau-Yosida
proximal mapping  of $\theta$, and ${\rm Pr}_{\theta^*}:\Re^{m\times n}\to\Re^{m\times n}$ is the Moreau-Yosida
proximal mapping of $\theta^*$. Denote $X:=\overline{X}+\overline{S}$. Let $X$ admit the following singular value decomposition 
\begin{equation}\label{eq:SVD-X}
X=\overline{U}\left[ \Sigma(X)\ \  0 \right]\overline{V}^{T}.
\end{equation}
Let $\sigma=\sigma(X)$, $\overline{\sigma}=\sigma(\overline{X})$ and $\overline{u} = \sigma(\overline{S})$ be the singular values of $X$, $\overline{X}$ and $\overline{S}$, respectively. Since $\|\cdot\|$ and $\|\cdot\|_{(k)}$ are unitarily invariant, we know from von Neumann's trace inequality (Lemma \ref{lem:vonNeumann}) that
\begin{equation}\label{eq:MY}
\overline{X}=\overline{U}\left[{\rm Diag}\,(\overline{\sigma})\ \ 0\right]\overline{V}^{T} \quad {\rm and} \quad \overline{S}=\overline{U}\left[{\rm Diag}\,(\overline{u})\ \ 0\right]\overline{V}^{T} \quad \mbox{with $\overline{\sigma}=g(\sigma)$ and $\overline{u}=\sigma-g(\sigma)$}, 
\end{equation}
where $g:\Re^{m}\to\Re^m$ is the Moreau-Yosida proximal mapping of the vector $k$-norm (i.e., the sum of the $k$ largest components in absolute value of any vector in $\Re^{m}$). The properties of the proximal mapping $g$ have been studied recently in \cite{WDSToh11}, e.g., for any given $x\in\Re^{m}$, the unique optimal solution $g(x)\in\Re^{m}$ can be computed within $O(m)$ arithmetic operations (see \cite[Section 3.1]{WDSToh11} for details). The following simple observations are useful for our subsequence analysis, which can be obtained  directly from the characterization of the subdifferential of $\theta=\|\cdot\|_{(k)}$ (cf. \cite{Watson92,OWomersley93}).

\begin{lemma}\label{lem:MY-vector-KKT}
	$\overline{\sigma}$ and $\overline{u}$ are the singular values of the solution $(\overline{X},\overline{S})$ of the GE \eqref{eq:GE} (or \eqref{eq:GE-dual}) if and only if $\overline{\sigma}$ and $\overline{u}$ satisfy the following conditions. 
	\begin{itemize}
		\item[(i)] If $\overline{\sigma}_{k}>0$, then 
		\begin{equation}\label{eq:condition-u-1}
		\overline{u}_{\alpha}=e_{\alpha}, \quad 0\le \overline{u}_{\beta}\le e_{\beta}, \quad \sum_{i\in\beta}\overline{u}_{i}=k-k_{0} \quad {\rm and} \quad \overline{u}_{\gamma}=0,
		\end{equation}
		where $0\leq k_{0}\leq k-1$ and $k\leq
		k_{1}\leq m$ are two integers  such that
		\begin{equation}\label{eq:def-k1-k2-case1}
		\overline{\sigma}_{1}\ge \ldots\ge
		\overline{\sigma}_{k_{0}}>\overline{\sigma}_{k_{0}+1}=\ldots=\overline{\sigma}_{k}=\ldots=\overline{\sigma}_{k_{1}}>\overline{\sigma}_{k_{1}+1}\ge
		\ldots\ge \overline{\sigma}_{m}\ge 0
		\end{equation}
		and
		\begin{equation}\label{eq:def-alpha-beta-gamma-case1}
		\alpha=\{1,\ldots,k_{0} \},\quad \beta=\{k_{0}+1,\ldots,k_{1}\}\quad {\rm and} \quad \gamma=\{k_{1}+1,\ldots,m\}.
		\end{equation}
		
		\item[(ii)] If $\overline{\sigma}_{k}=0$, then 
		\begin{equation}\label{eq:condition-u-2}
		\overline{u}_{\alpha}=e_{\alpha}, \quad 0\le \overline{u}_{\beta}\le e_{\beta}\quad {\rm and}\quad  \sum_{i\in\beta}\overline{u}_{i}\leq k-k_{0},
		\end{equation}
		where $0\leq k_{0}\leq k-1$ is the integer such that
		\begin{equation}\label{eq:def-k1-k2-case2}
		\overline{\sigma}_{1}\ge \cdots \ge \overline{\sigma}_{k_{0}}>\overline{\sigma}_{k_{0}+1}=\ldots=\overline{\sigma}_{k}=\ldots=\overline{\sigma}_{m}=0
		\end{equation}
		and
		\begin{equation}\label{eq:def-alpha-beta-case2}
		\alpha=\{1,\ldots,k_{0} \}\quad {\rm and} \quad \beta=\{k_{0}+1,\ldots,m\}.
		\end{equation}
	\end{itemize}
	%	Moreover, if $\overline{X}\neq 0$,
	%	\begin{equation}\label{eq:dual-active}
	%	\|\overline{S}\|^*_{(k)}=\max\{\|\overline{S}\|_2,\|\overline{S}\|_*/k \}=1.
	%	\end{equation}
\end{lemma}

For notational convenience, we use $\beta_1$, $\beta_2$ and $\beta_3$ to denote the index sets 
\begin{equation}\label{eq:def-beta123-case1}
\beta_{1}:=\{i\in\beta\,|\, \overline{u}_{i}=1  \},\quad \beta_{2}:=\{i\in\beta\,|\, 0<\overline{u}_{i}<1  \}\quad {\rm and}\quad \beta_{3}:=\{i\in\beta\,|\, \overline{u}_{i}=0  \}.
\end{equation}
For $X=\overline{X}+\overline{S}$, let $a$, $b$ and $c$ be the index sets defined by \eqref{eq:def-a-b-c}. We use $a_1,\ldots,a_r$ to denote the index sets defined by \eqref{eq:ak-nonsymmetric} with respect to $X$ and $a_{r+1}=b$ for the sake of convenience. Thus, by Lemma \ref{lem:MY-vector-KKT}, we know that if $\overline{\sigma}_{k}>0$,  then there exist integers $r_{0} \le r_{1}\in\{0,1,\ldots,r+1\}$, $r_{0} \le\widetilde{r}_0\le r_{0}+1$ and $r_{1}-1 \le\widetilde{r}_1\le r_1$ such that
\begin{equation}\label{eq:widetilde_r-case1}
\alpha=\bigcup_{l=1}^{r_{0}}a_{l}, \quad \beta_1=\bigcup_{l=r_0+1}^{\widetilde{r}_0}a_{l}, \quad \beta_2=\bigcup_{l=\widetilde{r}_0+1}^{\widetilde{r}_1}a_{l}, \quad \beta_3=\bigcup_{l=\widetilde{r}_1+1}^{r_1}a_{l} \quad {\rm and} \quad \gamma=\bigcup_{l=r_{1}+1}^{r+1}a_{l};
\end{equation}
if $\overline{\sigma}_{k}=0$, then there exist integers $r_{0}\in\{0,1,\ldots,r+1\}$ and $r_0\le\widetilde{r}_0\le r_0+1$ such that
\begin{equation}\label{eq:widetilde_r-case2}
\alpha=\bigcup_{l=1}^{r_{0}}a_{l}, \quad \beta_1=\bigcup_{l=r_{0}+1}^{\widetilde{r}_0}a_{l}, \quad \beta_2=\bigcup_{l=\widetilde{r}_0+1}^{r}a_{l}\quad {\rm and} \quad  \beta_{3}=b.
\end{equation}
Moreover, we know that for each $l\in\{1,\ldots,r_0\}$, $\overline{\sigma}_i=\overline{\sigma}_j$ for any $i,j\in a_l$, which implies that we can use $\overline{\nu}_1>\ldots>\overline{\nu}_{r_0}>0$ to denote those common values. Similarly, if $\overline{\sigma}_k>0$, we use $\overline{\mu}_{\widetilde{r}_0+1}>\ldots>\overline{\mu}_{\widetilde{r}_1}>0$ to denote the corresponding common values of $\overline{u}$; if $\overline{\sigma}_k=0$, we use $\overline{\mu}_{\widetilde{r}_0+1}>\ldots>\overline{\mu}_{r}>0$ to denote the corresponding common values of $\overline{u}$.

\section{The nondegeneracy and strict complementarity}\label{sect:nondegeneracy}
In this section, we shall introduce the nondegeneracy and strict complementarity of the optimization problem \eqref{eq:P-g}. To do so, let us consider the following conic reformulation of \eqref{eq:P-g}: 
\begin{equation}\label{eq:CP}
\begin{array}{cl} \min &  f(x)+t
\\ [3pt]
{\rm s.t.} & (g(x),t)\in {\cal K},
\end{array}
\end{equation}
where ${\cal K}={\rm epi}\,\theta$. 

Let $(\bar{x},\bar{t})$ be a feasible point of \eqref{eq:CP}. Denote $\overline{X}=g(\bar{x})\in\Re^{m\times n}$. Recall the definition \cite[Definition 6.1]{RWets98} of the tangent cone ${\cal T}_{{\cal K}}(\overline{X},\bar{t})$ of ${\cal K}$ at the given point $(\overline{X},\bar{t})\in{\cal K}$, i.e.,
\[
{\cal T}_{{\cal K}}(\overline{X},\bar{t})=\left\{(H,\tau)\in\Re^{m\times n}\times\Re\,|\, \exists\, \rho_{n}\downarrow 0,\ {\rm dist}\left((\overline{X},\bar{t})+\rho_{n}(H,\tau), {\cal K}\right)=o(\rho_{n})  \right\}.
\]
For any convex function $\phi:\Re^{m\times n}\to(-\infty,\infty)$, we know from \cite[Theorem 2.4.9]{Clarke83} that  
\begin{equation}\label{eq:T-epi-f}
{\cal T}_{{\rm epi}\,\phi}(Y,\phi(Y))= {\rm epi}\,\phi'(Y;\cdot):= \Big\{ (H,\tau)\in\Re^{m\times n}\times\Re\,|\,
\phi'(Y;H)\le \tau\Big\}, \quad Y\in\Re^{m\times n}.
\end{equation}
Therefore, for $\theta=\|\cdot\|_{(k)}$, we know from Proposition \ref{prop:easy} that  
\begin{equation}\label{eq:dirc-diff-k-norm}
{\cal T}_{\cal K}(\overline{X},\theta(\overline{X}))=\left\{
\begin{array}{ll}
\big\{(H,\tau) \mid {\rm tr}(\overline{U}_{\alpha}^{T}H\overline{V}_{\alpha})+\displaystyle{\sum_{i=1}^{k-k_{0}}}\lambda_{i}\left(S(\overline{U}_{\beta}^{T}H\overline{V}_{\beta})\right)\le \tau \big\} & \mbox{if $\sigma_k(\overline{X})>0$,} \\ [3pt]
\big\{ (H,\tau) \mid {\rm tr}(\overline{U}_{\alpha}^{T}H\overline{V}_{\alpha})+\displaystyle{\sum_{i=1}^{k-k_{0}}}\sigma_{i}\left(\left[\overline{U}_{\beta}^{T}H\overline{V}_{\beta}\
\ \overline{U}_{\beta}^{T}H\overline{V}_{2}\right]\right)\le \tau \big\} & \mbox{if $\sigma_k(\overline{X})=0$.}
\end{array}
\right.
\end{equation}
Define $G:{\cal X}\times\Re\to\Re^{m\times n}\times\Re$ by $G(x,t):=(g(x),t)$, $(x,t)\in{\cal X}\times\Re$. Robinson's CQ \cite{Robinson76} for \eqref{eq:CP} at a given feasible point $(\bar{x},\bar{t})$  can be written as
\begin{equation}\label{eq:def-Robinson-CQ-MCP}
G'(\bar{x},\bar{t})({\cal X}\times \Re)+ {\cal T}_{\cal K}(\overline{X},\bar{t})= \Re^{m\times n}\times\Re.
\end{equation}

\begin{proposition}\label{prop:RCQ}
	For any $\bar{x}\in{\cal X}$, Robinson's CQ \eqref{eq:def-Robinson-CQ-MCP} for \eqref{eq:CP} holds at $(\bar{x},\theta(g(\bar{x})))$. 
\end{proposition}
\noindent {\bf Proof.} Note that the directional derivative $\theta'(\overline{X};\cdot)$ of the Ky Fan $k$-norm is finite everywhere. Therefore, the results can be derived directly from \eqref{eq:def-Robinson-CQ-MCP} and \eqref{eq:T-epi-f}. In fact, we only need to show that for any given $(X,t)\in\Re^{m\times n}\times\Re$, there exists $(h,\eta)\in{\cal X}\times\Re$ and $(H,\tau)\in{\cal T}_{\cal K}(\overline{X},\bar{t})$ with $\bar{t}=\theta(g(\bar{x}))$ such that 
\[
(g'(\bar{x})h,\eta)+(H,\tau)=(X,t).
\]
Let $H=X$ and $\tau=\theta'(\overline{X};X)$. By choosing $h=0$ and $\eta=t-\tau$, we know that the above equality holds trivially.   $\hfill \Box$

\vskip 10 true pt

As we mentioned in Section \ref{sect:Intro}, for a locally optimal solution $\bar{x}$ to the optimization problem \eqref{eq:P-g}, the corresponding Lagrange multiplier always exists. In next proposition, we show  that the set of Lagrange multipliers of \eqref{eq:P-g} is also convex, bounded and compact.

\begin{proposition}
	Let $\bar{x}\in{\cal X}$ be a locally optimal solution to the problem \eqref{eq:P-g}. The set of Lagrange multipliers of \eqref{eq:P-g} is a nonempty, convex, bounded and compact subset of $\Re^{m\times n}$.
\end{proposition}
\noindent {\bf Proof.} It is easy to see that $\bar{x}\in{\cal X}$ is a locally optimal solution of \eqref{eq:P-g} if and only if $(\bar{x},\theta(g(\bar{x})))$ is a locally optimal solution of \eqref{eq:CP}. Moreover, by \eqref{eq:GE=CP}, we know that there exists a Lagrange multiplier $\overline{S}\in\Re^{m\times n}$ if and only if  there exists $\overline{S}\in\Re^{m\times n}$ such that the following KKT condition of \eqref{eq:CP} holds at $(\bar{x},\theta(g(\bar{x})),\overline{S},-1)$:
\begin{equation}\label{eq:KKT-CP}
\left\{
\begin{array}{l}
\nabla f(x)+g'(x)^*S = 0, \\ [3pt]
\xi +1 = 0, \\ [3pt]
\left\langle (g(x),t),(S,\xi)\right\rangle=0,\quad (g(x),t)\in{\cal K} \quad {\rm and} \quad  (S,\xi)\in {\cal K}^{\circ}. 
\end{array}
\right.
\end{equation}

On the other hand, it is well-known \cite{ZKurcyusz79} that for a locally optimal solution of \eqref{eq:CP}, the corresponding  set of Lagrange multipliers   is nonempty, convex, bounded and compact if and only if Robinson's CQ holds. Therefore, the result follows from Proposition \ref{prop:RCQ} directly. $\hfill \Box$

\vskip 10 true  pt
Next, let us study the concept of nondegeneracy for the optimization problem \eqref{eq:P-g}.
For any convex function $\phi:\Re^{m\times n}\to(-\infty,\infty)$ and $Y\in\Re^{m\times n}$, the lineality space of ${\cal T}_{{\rm epi}\,\phi}(Y,\phi(Y))$, i.e., the largest linear subspace in ${\cal T}_{{\rm epi}\,\phi}(Y,\phi(Y))$, can be written as 
\begin{eqnarray}
&&{\rm lin}\left({\cal T}_{{\rm epi}\,\phi}Y,\phi(Y))\right)={\cal T}_{{\rm epi}\,\phi}(Y,\phi(Y))\cap(-{\cal T}_{{\rm epi}\,\phi}(Y,\phi(Y))) \nonumber \\ [3pt]
&=& \Big\{ (H,\tau)\in\Re^{m\times n}\times\Re\,|\,
\phi'(Y;H)\leq \tau \leq -\phi'(Y;-H) \Big\} \nonumber \\ [3pt]
&=& \Big\{ (H,\tau)\in\Re^{m\times n}\times\Re\,|\,
\phi'(Y;H)=-\phi'(Y;-H)=\tau\Big\} \label{eq:lin-Tangent}.
\end{eqnarray}
The last equation of \eqref{eq:lin-Tangent} follows from \cite[Theorem 23.1]{Rockafellar70}, directly. For the Ky Fan $k$-norm $\theta=\|\cdot\|_{(k)}$, define the linear subspace  ${\cal T}^{\rm lin}(\overline{X})\subseteq\Re^{m\times n}$ by
\begin{equation}\label{eq:def-Tlin-A}
{\cal T}^{\rm lin}(\overline{X}):=\Big\{ 
H\in\Re^{m\times n}\,|\, \theta'(\overline{X};H)=-\theta'(\overline{X};-H)\Big\}.
\end{equation}
If $\overline{S}\in\partial\,\theta(\overline{X})$,
then, by Proposition \ref{prop:easy}, we have
\begin{equation}\label{eq:Tlin}
{\cal T}^{\rm lin}(\overline{X})=\left\{ 
\begin{array}{ll}
\left\{ H\in\Re^{m\times n}\mid  S(\overline{U}_{\beta}^{T}H\overline{V}_{\beta})= \tau I_{|\beta|}\ \mbox{for some $\tau\in\Re$} \right\}  & \mbox{if $\sigma_k(\overline{X})>0$,} \\ [5pt]
\left\{ H\in\Re^{m\times n}\mid \big[\overline{U}_{\beta}^{T}H\overline{V}_{\beta}\
\ \overline{U}_{\beta}^{T}H\overline{V}_{2}\big]=0\right\} & \mbox{if $\sigma_k(\overline{X})=0$,}
\end{array}\right.
\end{equation}
where $\overline{U}\in{\cal O}^m$ and $\overline{V}\in{\cal O}^n$ are eigenvectors of $X=\overline{X}+\overline{S}$, and the index set $\beta$ is defined in \eqref{eq:def-alpha-beta-gamma-case1} if $\sigma_k(\overline{X})>0$ and in \eqref{eq:def-alpha-beta-case2} if $\sigma_k(\overline{X})=0$.

For the problem \eqref{eq:CP}, the concept of Robinson's constraint nondegeneracy \cite{Robinson84,Robinson87} can be specified as follows. The constraint nondegeneracy for \eqref{eq:CP} holds at the feasible point $(\bar{x},\bar{t})$  if 
\begin{equation}\label{eq:def-constraint-nondegenerate-QMCP}
G'(\bar{x},\bar{t})({\cal X}\times\Re)+{\rm lin}\left({\cal T}_{\cal K}(\overline{X},\bar{t})\right)=\Re^{m\times n}\times\Re,
\end{equation}
where the lineality space ${\rm lin}\left({\cal T}_{\cal K}(\overline{X},\bar{t})\right)$ is given by \eqref{eq:lin-Tangent} with respect to $\theta=\|\cdot\|_{(k)}$.

\begin{proposition}\label{prop:RCN-primal}
	The constraint nondegeneracy \eqref{eq:def-constraint-nondegenerate-QMCP} for \eqref{eq:CP} holds at $(\bar{x},\theta(\overline{X}))$ if and only if
	\begin{equation}\label{eq:def-constraint-nondegenerate-QMOP}
	g'(\bar{x}){\cal X} + {\cal T}^{\rm lin}(\overline{X}) = \Re^{m\times n},
	\end{equation}
	where ${\cal T}^{\rm lin}(\overline{X})\in\Re^{m\times n}$ is the linear subspace defined by \eqref{eq:def-Tlin-A}. Therefore, we say that the nondegeratacy for the problem \eqref{eq:P-g} holds at $\bar{x}$  if   \eqref{eq:def-constraint-nondegenerate-QMOP} holds. 
\end{proposition}
\noindent {\bf Proof.} For any given $X\in\Re^{m\times n}$, by \eqref{eq:def-constraint-nondegenerate-QMCP}, we know that there exists $h\in{\cal X}$, $(H,\eta)\in{\rm lin}\left({\cal T}_{\cal K}(\overline{X},\bar{t})\right)$ such that
\[
(g'(\bar{x})h,-\eta)+(H,\eta) = (X,0).
\]
Since $H\in{\cal T}^{\rm lin}(\overline{X})$, we know that \eqref{eq:def-constraint-nondegenerate-QMOP} holds. 

Conversely, for any $(X,t)\in\Re^{m\times n}\times\Re$, by \eqref{eq:def-constraint-nondegenerate-QMOP}, we know that there exists $h\in{\cal X}$ and $H\in{\cal T}^{\rm lin}(\overline{X})$ such that 
\[
g'(\bar{x})h+H=X.
\]
Denote $\tau=\theta'(\overline{X};H)$. By taking $\eta=t-\tau$, we obtain that
\[
(g'(\bar{x})h,\eta)+(H,\tau) = (X,t),
\]
which implies that the constraint nondegeneracy \eqref{eq:def-constraint-nondegenerate-QMCP} holds at $(\bar{x},\theta(\overline{X}))$.  $\hfill \Box$

\vskip  10 true pt
Let $\bar{x}\in{\cal X}$ be a locally optimal solution of \eqref{eq:P-g}. Denote $\overline{X}=g(\bar{x})$. Let $\beta$ be the index set defined in \eqref{eq:def-alpha-beta-gamma-case1} if $\sigma_k(\overline{X})>0$ and in \eqref{eq:def-alpha-beta-case2} if $\sigma_k(\overline{X})=0$.  The following definition of the strict complementarity of \eqref{eq:P-g} can be regarded as a generalization of the strict complementarity for the constraint optimization problem (cf. \cite[Definition 4.74]{BShapiro00}).
\begin{definition}
	We say the strict complementarity condition holds at $\bar{x}\in{\cal X}$ if there exists $\overline{S}\in{\rm ri}\,(\partial\,\theta(\overline{X}))$ such that 
	\begin{equation}\label{eq:KKT-eq}
	\nabla f(\bar{x})+g'(\bar{x})^*\overline{S}=0.
	\end{equation}
\end{definition}

By Lemma \ref{lem:MY-vector-KKT}, one can derive the following proposition easily. For simplicity, we omit the detail proof here.
\begin{proposition}\label{prop:strict complementarity}
	The strict complementarity condition holds at $\bar{x}\in{\cal X}$ if and only if there exists $\overline{S}\in\partial\,\theta(\overline{X})$ such that \eqref{eq:KKT-eq} holds and
	\begin{itemize}
		\item[(i)] if $\sigma_k(\overline{X})>0$, then $0<\sigma_{\beta}(\overline{S})<e_{\beta}$;
		\item[(ii)] if $\sigma_k(\overline{X})=0$, then $\sigma_{\beta}(\overline{S})<e_{\beta}$ and $\sum_{i\in\beta}\sigma_i(\overline{S})<k-k_0$,
	\end{itemize} 
\end{proposition}

\begin{proposition}
	Let $\bar{x}\in{\cal X}$ be a locally optimal solution of \eqref{eq:P-g}. Denote $\overline{X}=g(\bar{x})$. If $\bar{x}$ is nondegenerate, then $\overline{S}$ satisfying \eqref{eq:KKT-g} is unique. Conversely, if $\overline{S}$ satisfying \eqref{eq:KKT-g} is unique and the strict complementarity condition holds at $\bar{x}$, then $\bar{x}$ is nondegenerate.
\end{proposition}
\noindent {\bf Proof.}  The following proof is a slight modification of the proof of \cite[Proposition 4.75]{BShapiro00}. Suppose that $\bar{x}$ is nondegenerate and let $\overline{S}$ and $\overline{S}'$ satisfy \eqref{eq:KKT-g}. Then, we know that $g'(\bar{x})^*(\overline{S}-\overline{S}')=0$, which implies that $\Delta:=\overline{S}-\overline{S}'\in \left[g'(\bar{x}){\cal X}\right]^{\perp}$. Denote $X=\overline{X}+\overline{S}$ and $X'=\overline{X}+\overline{S}'$. Suppose that $X$ and $X'$ admit the SVD:
\[
X=\overline{U}[\Sigma(X)\ \ 0]\overline{V}^T\quad {\rm and} \quad X'=\overline{U}'[\Sigma(X')\ \ 0](\overline{V}')^T,
\]
where $\overline{U},\overline{U}'\in{\cal O}^m$ and $\overline{V},\overline{V}'\in{\cal O}^n$. By \eqref{eq:MY}, we know that both $(\overline{U},\overline{V})$ and $(\overline{U}',\overline{V}')$ are eigenvalue vectors of $\overline{X}$. Therefore, it follows from \cite[Proposition 5]{DSToh10} that if $\sigma_k(\overline{X})>0$, then there exist orthogonal matrices $Q_1\in{\cal O}^{|\alpha|}$, $Q_2\in{\cal O}^{|\beta|}$, $Q_3\in{\cal O}^{|\gamma|}$ and ${\cal Q}_3'\in{\cal O}^{|\gamma|+n-m}$ such that 
\[
\overline{U}'=\overline{U}\left[\begin{array}{ccc}
Q_1 & 0 & 0 \\
0 & Q_2 & 0 \\
0 & 0 & Q_3
\end{array}\right] \quad {\rm and} \quad \overline{V}'=\overline{V}\left[\begin{array}{ccc}
Q_1 & 0 & 0 \\
0 & Q_2 & 0 \\
0 & 0 & Q_3'
\end{array}\right];
\]
if $\sigma_k(\overline{X})=0$, then there exist orthogonal matrices $Q_1\in{\cal O}^{|\alpha|}$, $Q_2\in{\cal O}^{|\beta|}$ and $Q_2'\in{\cal O}^{|\beta|+n-m}$ such that 
\[
\overline{U}'=\overline{U}\left[\begin{array}{cc}
Q_1 & 0 \\
0 & Q_2 
\end{array}\right] \quad {\rm and} \quad \overline{V}'=\overline{V}\left[\begin{array}{ccc}
Q_1 & 0 \\
0 & Q_2'
\end{array}\right].
\]  
Therefore, by \eqref{eq:MY}, we know from Lemma \ref{lem:MY-vector-KKT} that if $\sigma_k(\overline{X})>0$, then
\begin{equation}\label{eq:H1-tilde-case1-dual}
\overline{U}^TH\overline{V} = \left[
\begin{array}{cccc}
0 & 0 & 0 & 0 \\ [3pt]
0 & S(\overline{U}^T_{\beta}H\overline{V}_{\beta}) & 0 & 0 \\ [3pt]
0 & 0 & 0 & 0 
\end{array}
\right] \quad {\rm with} \quad {\rm tr}(S(\overline{U}^T_{\beta}H\overline{V}_{\beta}))=0.
\end{equation}
if $\sigma_k(\overline{X})=0$, then
\begin{equation}\label{eq:H1-tilde-case2-dual}
\overline{U}^TH\overline{V} = \left[
\begin{array}{ccc}
0 & 0 & 0 \\ [3pt]
0 & \overline{U}^T_{\beta}H\overline{V}_{\beta} & \overline{U}^T_{\beta}H\overline{V}_{2}
\end{array}
\right].
\end{equation}
Thus, we know from \eqref{eq:Tlin} that in both cases,
\[
\langle \Delta,H\rangle = \langle \overline{U}^T\Delta \overline{V},\overline{U}^TH\overline{V}\rangle = 0 \quad \forall\, H\in {\cal T}^{\rm lin}(\overline{X}),
\]
which implies that $\Delta\in \left[{\cal T}^{\rm lin}(\overline{X})\right]^{\perp}$. Therefore, by \eqref{eq:def-constraint-nondegenerate-QMOP}, we know that $\Delta =0$, i.e., $\overline{S}$ satisfying \eqref{eq:KKT-g} is unique.

Conversely, since the strict complementarity condition holds at $\bar{x}$, we know that the unique Lagrange multiplier $\overline{S}\in\partial\,\theta(\overline{X})$ satisfying (i) and (ii) of Proposition \ref{prop:strict complementarity}. Let $X=\overline{X}+\overline{S}$ admit the SVD \eqref{eq:SVD-X}. Suppose that the constraint nondegenerate condition \eqref{eq:def-constraint-nondegenerate-QMOP} does not hold at $\overline{X}$, i.e., there exists $0\neq H\in\left[g'(\bar{x}){\cal X}\right]^{\perp}\cap\left[{\cal T}^{\rm lin}(\overline{X})\right]^{\perp}$. Therefore, we know that $g'(\overline{X})^*H=0$. Moreover, by \eqref{eq:Tlin}, we know that if $\sigma_k(\overline{X})>0$, then \eqref{eq:H1-tilde-case1-dual} holds; if $\sigma_k(\overline{X})=0$, then \eqref{eq:H1-tilde-case2-dual} holds. Since $g'(\overline{X})^*H=0$, we know that for any $\rho$,
\[
\nabla f(\bar{x})+g'(\overline{X})^*(\overline{S}+\rho H)=0.
\]
Moreover, since $\overline{S}$ satisfies (i) and (ii) of Proposition \ref{prop:strict complementarity}, by \eqref{eq:H1-tilde-case1-dual} and \eqref{eq:H1-tilde-case2-dual}, we know from Lemma \ref{lem:MY-vector-KKT} that for $\rho>0$ small enough, $\overline{S}+\rho H\in\partial\,\theta(\overline{X})$. This contradicts the uniqueness of $\overline{S}$. $\hfill \Box$

\begin{remark}
	Let $\overline{X}\in\partial\,\theta^*(\overline{S})$. For the dual norm $\vartheta=\|\cdot\|_{(k)}^*$, since $(\overline{S},-1)\in{\cal K}^{\circ}$, we have 
	\[
	{\cal T}_{{\cal K}^{\circ}}(\overline{S},-1)=\left\{
	\begin{array}{ll}
	\left\{(H,\tau)\in\Re^{m\times n}\times\Re\mid \vartheta'(\overline{S};H)\le -\tau\right\} & \mbox{if $\vartheta(\overline{S})=1$,} \\ [3pt]
	\Re^{m\times n}\times \Re & \mbox{if $\vartheta(\overline{S})<1$.}
	\end{array}
	\right.
	\]
	We define the linear subspace ${\cal T}_{\circ}^{\rm lin}(\overline{S})\subseteq\Re^{m\times n}$ by \begin{equation}\label{eq:def-Tlin-dual-A}
	{\cal T}_{\circ}^{\rm lin}(\overline{S}):=\left\{
	\begin{array}{ll}
	\big\{ 
	H\in\Re^{m\times n}\,|\, \vartheta'(\overline{S};H)=-\vartheta'(\overline{S};-H)=0\big\}  & \mbox{if $\vartheta(\overline{S})=1$,} \\ [3pt]
	\Re^{m\times n} & \mbox{if $\vartheta(\overline{S})<1$.}
	\end{array}
	\right.
	\end{equation}  
	For the case that $\vartheta(\overline{S})=\max\{\|\overline{S}\|_2,\|\overline{S}\|_*/k \}=1$, we know from Proposition \ref{prop:easy} that if $\|\overline{S}\|_*<k$, then
	\begin{equation}\label{eq:cha-Tlin-dual-1}
	{\cal T}_{\circ}^{\rm lin}(\overline{S})=\left\{ H\in\Re^{m\times n}\mid S([\overline{U}_{\alpha}\ \ \overline{U}_{\beta_1}]^TH[\overline{V}_{\alpha}\ \ \overline{V}_{\beta_1}])=0 \right\};
	\end{equation}
	if $\|\overline{S}\|_*=k$, then
	\begin{equation}\label{eq:cha-Tlin-dual-2} 
	{\cal T}_{\circ}^{\rm lin}(\overline{S})=\left\{H\in\Re^{m\times n}\mid S(\overline{U}_{\alpha\cup\beta_1}^TH\overline{V}_{\alpha\cup\beta_1})=0,\, {\rm
		tr}(\overline{U}_{\beta_2}^{T}H\overline{V}_{\beta_2})=0,\,  \overline{U}_{\beta_3\cup\gamma}^TH\overline{V}_{\beta_3\cup\gamma\cup c}=0\right\},
	\end{equation}
	where $\overline{U}\in{\cal O}^m$ and $\overline{V}\in{\cal O}^n$ are eigenvectors of $X=\overline{X}+\overline{S}$, the index set $\beta$ is defined in \eqref{eq:def-alpha-beta-gamma-case1} if $\sigma_k(\overline{X})>0$ and in \eqref{eq:def-alpha-beta-case2} if $\sigma_k(\overline{X})=0$, and $\beta_1$, $\beta_2$ and $\beta_3$ are the index sets defined by \eqref{eq:def-beta123-case1}.
\end{remark}

\section{The critical cones}\label{sect:critical cones}
From now on, let us always assume that $\overline{X}=g(\bar{x})$ and $\overline{S}$ are solutions of the GEs \eqref{eq:GE} and \eqref{eq:GE-dual}. Therefore, the critical cones associated with the GEs \eqref{eq:GE} and \eqref{eq:GE-dual} can be defined correspondingly from the critical cones  associated the complementarity problem \eqref{eq:GE=CP}.

Firstly, consider the GE \eqref{eq:GE}. Denote $(X,t)=(\overline{X}+\overline{S},\theta(\overline{X})-1)$. The critical cone of ${\cal K}$ at $(X,t)$ associated with the complementarity problem in \eqref{eq:GE=CP}, is defined as
\begin{equation}\label{eq:def-critical cone-K}
{\cal C}\left((X,t);{\cal K}\right)={\cal T}_{\cal K}( \overline{X},\theta(\overline{X}))\cap(\overline{S},-1)^{\perp}.
\end{equation}
Thus, we know from \eqref{eq:T-epi-f} that 
\begin{equation}\label{eq:eqiv-critical-1}
(H,\tau)\in{\cal C}\left((X,t);{\cal K}\right) \Longleftrightarrow 
\left\{
\begin{array}{l}
H\in{\cal C}(X;\partial\,\theta(\overline{X})), \\ [3pt]
\tau = \langle \overline{S},H\rangle,
\end{array}
\right. 
\end{equation}
where ${\cal C}(X;\partial\,\theta(\overline{X}))\subseteq\Re^{m\times n}$ is defined by
\begin{equation}\label{eq:def-critical-cone}
{\cal C}(X;\partial\,\theta(\overline{X})):=\left\{H\in\Re^{m\times n}\mid \theta'(\overline{X};H)\le \langle \overline{S},H\rangle \right\}.
\end{equation}
Since $\theta'(\overline{X};\cdot)$ is a positively homogeneous convex function with $\theta'(\overline{X};0)=0$, ${\cal C}(X;\partial\,\theta(\overline{X}))$ is indeed a closed convex cone. We call ${\cal C}(X;\partial\,\theta(\overline{X}))$ the critical cone of $\partial\,\theta(\overline{X})$ at $X=\overline{X}+\overline{S}$, associated with the GE \eqref{eq:GE}.

Next, we present the following proposition on the characterization of the critical cone ${\cal C}(X;\partial\,\theta(\overline{X}))$. 

\begin{proposition}\label{prop:critical-cone}
	Suppose that $(\overline{X},\overline{S})\in\Re^{m\times n}\times \Re^{m\times n}$ is a solution of the GE \eqref{eq:GE}. Let $X=\overline{X}+\overline{S}$ admit the SVD \eqref{eq:SVD-X}. Then,
	\begin{equation}\label{eq:eqiv-critical-2}
	H\in{\cal C}(X;\partial\,\theta(\overline{X}))\quad \Longleftrightarrow \quad  \theta'(\overline{X};H)=\langle \overline{S},H\rangle,
	\end{equation} 
	which is equivalent to the following conditions. 
	\begin{itemize}
		\item[(i)] If $\sigma_k(\overline{X})>0$, then there exists some $\tau\in\Re$ such that 
		\[ 
		\lambda_{|\beta_1|}(S(\overline{U}_{\beta_1}^TH\overline{V}_{\beta_1}))\ge \tau \ge \lambda_1(S(\overline{U}_{\beta_3}^TH\overline{V}_{\beta_3})) 
		\] and
		\[
		S(\overline{U}_{\beta}^TH\overline{V}_{\beta})= \left[
		\begin{array}{ccc}
		S(\overline{U}_{\beta_1}^TH\overline{V}_{\beta_1}) & 0 & 0\\
		0 & \tau I_{|\beta_2|} & 0  \\
		0 & 0 & S(\overline{U}_{\beta_3}^TH\overline{V}_{\beta_3}) 
		\end{array}
		\right].
		\]
		\item[(ii)] If $\sigma_k(\overline{X})=0$ and $\|\overline{S}\|_*=k$, then there exists some $\tau\ge 0$ such that 
		\[
		\lambda_{|\beta_1|}(S(\overline{U}_{\beta_1}^TH\overline{V}_{\beta_1}))\ge \tau \ge \sigma_1\left(\big[\overline{U}_{b}^{T}H\overline{V}_{b}\
		\ \overline{U}_{b}^{T}H\overline{V}_{2}\big]\right) 
		\] and
		\[
		\left[\overline{U}_{\beta}^{T}H\overline{V}_{\beta}\ \ \overline{U}_{\beta}^{T}H\overline{V}_{2}\right]=\left[\begin{array}{cccc}
		S(\overline{U}_{\beta_1}^{T}H\overline{V}_{\beta_1}) & 0 & 0 & 0 \\
		0 & \tau I_{|\beta_2|} &  0 & 0 \\
		0  & 0 &  \overline{U}_{b}^{T}H\overline{V}_{b} & \overline{U}_{b}^{T}H\overline{V}_{2} \end{array}\right].
		\]
		\item[(iii)] If $\sigma_k(\overline{X})=0$ and $\|\overline{S}\|_*<k$, then $S(\overline{U}_{\beta_1}^{T}H\overline{V}_{\beta_1})\succeq 0$ and
		\[
		\left[\overline{U}_{\beta}^{T}H\overline{V}_{\beta}\ \ \overline{U}_{\beta}^{T}H\overline{V}_{2}\right]=\left[\begin{array}{cccc}
		S(\overline{U}_{\beta_1}^{T}H\overline{V}_{\beta_1}) & 0 & 0 & 0 \\
		0 & 0 &  0 & 0 \\
		0  & 0 &  0 & 0 \end{array}\right].
		\]
	\end{itemize}
	
\end{proposition}
\noindent{\bf Proof.} Denote $\overline{\sigma}=\sigma(\overline{X})$ and $\overline{u}=\sigma(\overline{S})$. By
\[
\left\langle \overline{S},H\right\rangle=\left\langle \overline{U}^T\overline{S}\overline{V},\overline{U}^TH\overline{V}\right\rangle = \left\langle [{\rm Diag}(\overline{u})\ \ 0],\overline{U}^TH\overline{V}\right\rangle,
\] 
we know from Lemma \ref{lem:MY-vector-KKT} that for any $H\in\Re^{m\times n}$,
\[
\left\langle \overline{S},H \right\rangle = \left\{
\begin{array}{ll}
 {\rm tr}( \overline{U}_{\alpha}^{T}H\overline{V}_{\alpha})+\left\langle  {\rm Diag}(\overline{u}_{\beta}), S(\overline{U}_{\beta}^{T}H\overline{V}_{\beta})\right\rangle & \mbox{if $\overline{\sigma}_k>0$,} \\ [5pt]
 {\rm tr}( \overline{U}_{\alpha}^{T}H\overline{V}_{\alpha})+\left\langle [{\rm Diag}(\overline{u}_{\beta})\ \ 0], \left[\overline{U}_{\beta}^{T}H\overline{V}_{\beta}\
 \ \overline{U}_{\beta}^{T}H\overline{V}_{2}\right]\right\rangle & \mbox{if $\overline{\sigma}_k=0$.}
\end{array}
\right.
\]
Thus, by combining with Fan's inequality (Lemma \ref{lem:Fan}) and  von Neumann's trace inequality (Lemma \ref{lem:vonNeumann}) , we obtain that for any $H\in\Re^{m\times n}$, if $\overline{\sigma}_k>0$,
\begin{equation}\label{eq:critical-cone-eq-1}
\left\langle  {\rm Diag}(\overline{u}_{\beta}), S(\widetilde{H}_{\beta\beta})\right\rangle \le \overline{u}_{\beta}^T\lambda(S(\widetilde{H}_{\beta\beta}))\le \displaystyle{\sum_{i=1}^{k-k_{0}}}\lambda_{i}\left(S(\widetilde{H}_{\beta\beta})\right),
\end{equation} 
and if $\overline{\sigma}_k=0$,
\begin{equation}\label{eq:critical-cone-eq-2}
\left\langle \big[{\rm Diag}(\overline{u}_{\beta})\ \ 0\big], \big[\widetilde{H}_{\beta\beta}\
\ \widetilde{H}_{\beta c}\big]\right\rangle  \le  \overline{u}_{\beta}^T\sigma\left(\big[\widetilde{H}_{\beta\beta}\
\ \widetilde{H}_{\beta c}\big]\right) \le  \displaystyle{\sum_{i=1}^{k-k_{0}}}\sigma_{i}\left(\big[\widetilde{H}_{\beta\beta}\
\ \widetilde{H}_{\beta c}\big]\right),
\end{equation}
where $\widetilde{H}=\overline{U}^TH\overline{V}$.
Therefore, we know from \eqref{eq:dirc-diff-k-norm}  that
\begin{equation*}
H\in{\cal C}(X;\partial\,\theta(\overline{X})) \Longleftrightarrow  \theta'(\overline{X};H)=\langle \overline{S},H\rangle \Longleftrightarrow \mbox{the equalities in \eqref{eq:critical-cone-eq-1} and \eqref{eq:critical-cone-eq-2} hold.}
\end{equation*}
Consider the following two cases.

{\bf Case 1} $\overline{\sigma}_k>0$. It follows from Lemma \ref{lem:Fan} that the first equality of \eqref{eq:critical-cone-eq-1} holds if and only if ${\rm Diag}(\overline{u}_{\beta})$ and $S(\widetilde{H}_{\beta\beta})$ admit a simultaneous ordered eigenvalue decomposition, i.e.,  there exists $R\in{\cal O}^{|\beta|}$ such that
\begin{equation}\label{eq:simultaneous-ordered}
{\rm Diag}(\overline{u}_{\beta}) = R {\rm Diag}(\overline{u}_{\beta}) R^T \quad {\rm and} \quad S(\widetilde{H}_{\beta\beta}) = R \Lambda(S(\widetilde{H}_{\beta\beta}))R^T.
\end{equation}
Let $r_{0} \le r_{1}\in\{0,1,\ldots,r+1\}$, $r_{0} \le\widetilde{r}_0\le r_{0}+1$ and $r_{1}-1 \le\widetilde{r}_1\le r_1$ be the integers  such that \eqref{eq:widetilde_r-case1} holds.  Therefore, the orthogonal matrix $R\in{\cal O}^{|\beta|}$ has the following block diagonal structure:
\begin{equation}\label{eq:R-block-diag-1}
R = \left[
\begin{array}{ccc}
R_{1} & 0 & 0\\
0 & R_{2} & 0  \\
0 & 0 & R_{3} 
\end{array}
\right] \quad {\rm with} \quad 
R_{2}=\left[
\begin{array}{cccc}
R_{2}^{(1)} & 0 &0 \\
0 & \ddots & 0\\
0 & 0 & R_{2}^{(\widetilde{r}_1-\widetilde{r}_0)}
\end{array}
\right] ,
\end{equation}
where $R_{1}\in{\cal O}^{|\beta_1|}$, $R_{2}\in{\cal O}^{|\beta_2|}$, $R_3\in{\cal O}^{|\beta_3|}$ and $R_{2}^{(l)}\in{\cal O}^{|a_{\widetilde{r}_0+l}|}$, $l=1,\ldots,\widetilde{r}_1-\widetilde{r}_0$. Thus,  \eqref{eq:simultaneous-ordered} holds if and only if  $S(\widetilde{H}_{\beta\beta})\in{\cal S}^{|\beta|}$ has the following block diagonal structure:
\[
S(\widetilde{H}_{\beta\beta}) = \left[
\begin{array}{ccccc}
S(\widetilde{H}_{\beta_1\beta_1}) & 0 & \cdots & 0 & 0\\ [3pt]
0 & S(\widetilde{H}_{a_{\widetilde{r}_0+1}a_{\widetilde{r}_0+1}}) & \cdots & 0  & 0  \\ [3pt]
\vdots & \vdots  & \ddots & \vdots & \vdots \\ [3pt]
0 & 0 & \cdots & S(\widetilde{H}_{a_{\widetilde{r}_1}a_{\widetilde{r}_1}}) & 0 \\ [3pt]
0 & 0 & \cdots & 0  & S(\widetilde{H}_{\beta_3\beta_3})
\end{array}
\right] ,
\]
and the elements of $\left( \lambda(S(\widetilde{H}_{\beta_1\beta_1})),\lambda(S(\widetilde{H}_{a_{\widetilde{r}_0+1}a_{\widetilde{r}_0+1}})),\ldots, \lambda(S(\widetilde{H}_{a_{\widetilde{r}_1}a_{\widetilde{r}_1}})), \lambda(S(\widetilde{H}_{\beta_3\beta_3})) \right)$ are in non-increasing order and are the eigenvalues of the symmetric matrix $S(\widetilde{H}_{\beta\beta})$.

On the other hand, by \eqref{eq:condition-u-1}, we know that $\overline{u}_{\beta_1}=e_{\beta_1}$, $0< \overline{u}_{\beta_2}< e_{\beta_2}$, $\overline{u}_{\beta_3}=0$ and $\langle e_{\beta},\overline{u}_{\beta}\rangle=k-k_0$. Then, we can verify that the second equality of \eqref{eq:critical-cone-eq-1} holds if and only if 
\begin{equation}\label{eq:beta2-diag-cri}
\lambda_{i}(S(\widetilde{H}_{\beta\beta}))=\lambda_{j}(S(\widetilde{H}_{\beta\beta}))\quad \forall\, i,j\in \{|\beta_1|+1,\ldots,|\beta_1|+|\beta_2| \}.
\end{equation}
In fact, it is clear that \eqref{eq:beta2-diag-cri} implies the second equality of \eqref{eq:critical-cone-eq-1} holds. Conversely, without loss of generality, assume that $\beta_2\neq\emptyset$, then $k-k_0\in \{|\beta_1|+1,\ldots,|\beta_1|+|\beta_2| \}$. Suppose that there exists $i \in\left\{|\beta_1|+1,\ldots,|\beta_1|+|\beta_2| \right\}$ but $i\neq k-k_0$ such that   $\lambda_{i}(S(\widetilde{H}_{\beta\beta}))>\lambda_{k-k_0}(S(\widetilde{H}_{\beta\beta}))$ or $\lambda_{k-k_0}(S(\widetilde{H}_{\beta\beta}))>\lambda_{i}(S(\widetilde{H}_{\beta\beta}))$. Then, since $0< \overline{u}_{\beta_2}< e_{\beta_2}$ and $\langle e_{\beta},\overline{u}_{\beta}\rangle=k-k_0$, for both cases, we always have 
\begin{eqnarray*}
&&\sum_{i=1}^{k-k_0}\lambda_i(S(\widetilde{H}_{\beta\beta}))-\overline{u}_{\beta}^T\lambda(S(\widetilde{H}_{\beta\beta})) \\ [3pt] &=&\sum_{i=1}^{k-k_0}\lambda_i(S(\widetilde{H}_{\beta\beta}))(1-(\overline{u}_{\beta})_i)-\sum_{i=k-k_0+1}^{|\beta|}(\overline{u}_{\beta})_i\lambda_i(S(\widetilde{H}_{\beta\beta})) \\ [3pt]
&>& \sum_{i=1}^{k-k_0}\lambda_{k-k_0}(S(\widetilde{H}_{\beta\beta}))(1-(\overline{u}_{\beta})_i)-\sum_{i=k-k_0+1}^{|\beta|}(\overline{u}_{\beta})_i\lambda_{k-k_0}(S(\widetilde{H}_{\beta\beta})) \\
&=& \lambda_{k-k_0}(S(\widetilde{H}_{\beta\beta}))\left(k-k_0-\sum_{i=1}^{k-k_0}(\overline{u}_{\beta})_i-\sum_{i=k-k_0+1}^{|\beta|}(\overline{u}_{\beta})_i\right)=0,
\end{eqnarray*}
which implies that the second equality of \eqref{eq:critical-cone-eq-1} does not hold, which contradicts the assumption. Therefore, we know that $H\in{\cal C}(X;\partial\,\theta(\overline{X}))$ if and only if (i) holds.

{\bf Case 2} $\overline{\sigma}_k=0$. We know from Lemma \ref{lem:vonNeumann} that the first equality of \eqref{eq:critical-cone-eq-2} holds if and only if $[{\rm Diag}(\overline{u}_{\beta})\ \ 0]$ and $\big[\widetilde{H}_{\beta\beta}\
\ \widetilde{H}_{\beta c}\big]$ admit a simultaneous ordered SVD, i.e., there exist orthogonal matrices $E\in{\cal O}^{|\beta|}$ and $F\in{\cal O}^{|\beta|+n-m}$ such that
\begin{equation}\label{eq:simultaneous-ordered-2}
[{\rm Diag}(\overline{u}_{\beta})\ \ 0]=E[{\rm Diag}(\overline{u}_{\beta})\ \ 0]F^{T}\quad {\rm and}\quad  \big[\widetilde{H}_{\beta\beta}\
\ \widetilde{H}_{\beta c}\big]=E[\Sigma([\widetilde{H}_{\beta\beta}\
\ \widetilde{H}_{\beta c}])\ \ 0]F^{T}.
\end{equation}
Let $r_{0}\in\{0,1,\ldots,r+1\}$ and $r_0\le\widetilde{r}_0\le r_0+1$ be the integers such that \eqref{eq:widetilde_r-case2} holds.  Therefore, it follows from \cite[Proposition 5]{DSToh10} that there exist orthogonal matrices $Q_1\in{\cal O}^{|\beta_1|}$, $Q_2\in{\cal O}^{|\beta_2|}$, $Q_3\in{\cal Q}^{|\beta_3|}$ and $Q'_3\in{\cal O}^{|\beta_3|+n-m}$ such that
\begin{equation}\label{eq:R-block-diag-2}
E = \left[\begin{array}{ccc}
Q_1 & 0 & 0 \\ [3pt]
0 & Q_2 & 0 \\ [3pt]
0 & 0 & Q_3
\end{array}
\right] \quad {\rm and} \quad F = \left[\begin{array}{ccc}
Q_1 & 0 & 0 \\ [3pt]
0 & Q_2 & 0 \\ [3pt]
0 & 0 & Q'_3
\end{array}
\right] \quad {\rm with} \quad 
Q_{2}=\left[
\begin{array}{cccc}
Q_{2}^{(1)} & 0 &0 \\
0 & \ddots & 0\\
0 & 0 & Q_{2}^{(r-\widetilde{r}_0)}
\end{array}
\right] ,
\end{equation}
where $Q_{2}^{(l)}\in{\cal O}^{|a_{\widetilde{r}_0+l}|}$, $l=1,\ldots,r-\widetilde{r}_0$. Thus, \eqref{eq:simultaneous-ordered-2} holds if and only if $\big[\widetilde{H}_{\beta\beta}\
\ \widetilde{H}_{\beta c}\big]$ has the following block diagonal structure:
\begin{equation*} 
\big[\widetilde{H}_{\beta\beta}\ \ \widetilde{H}_{\beta c}\big]=\left[\begin{array}{cccccc}
\widetilde{H}_{a_{r_{0}+1}a_{r_{0}+1}} & 0 & \cdots & 0 & 0 & 0 \\
0 & \widetilde{H}_{a_{\widetilde{r}_{0}+1}a_{\widetilde{r}_{0}+1}} & \cdots & 0 & 0 & 0 \\
\vdots & \vdots & \ddots & \vdots & \vdots & \vdots \\
0 & 0 & \cdots & \widetilde{H}_{a_{r}a_{r}} & 0 & 0 \\
0 & 0 & \cdots & 0 &  \widetilde{H}_{bb} & \widetilde{H}_{bc} \end{array}\right]
\end{equation*} 
with  $\widetilde{H}_{a_{l}a_{l}}\in{\cal S}^{|a_l|}$, $l=r_{0}+1,\ldots,r$,  and the elements of 
\[
h:=\left( \lambda(\widetilde{H}_{a_{r_{0}+1}a_{r_{0}+1}}), \lambda(\widetilde{H}_{a_{\widetilde{r}_{0}+1}a_{\widetilde{r}_{0}+1}}),\ldots, \lambda(\widetilde{H}_{a_{r}a_{r}}), \sigma([\widetilde{H}_{bb}\ \ \widetilde{H}_{bc}]) \right)\in\Re^m
\] 
are nonnegative and in non-increasing order and   $h=\sigma\big([\widetilde{H}_{\beta\beta}\ \ \widetilde{H}_{\beta c}]\big)$.

On the other hand, by \eqref{eq:condition-u-2}, we know that $\overline{u}_{\beta_1}=e_{\beta_1}$, $0<\overline{u}_{\beta_2}<e_{\beta_2}$,  $\overline{u}_{\beta_3}=0$ and $\langle e_{\beta},\overline{u}_{\beta}\rangle\le k-k_0$. Then, we may conclude that the second equality of \eqref{eq:critical-cone-eq-2} holds if and only if 
\begin{equation}\label{eq:beta2-diag-cri-2}
\left\{\begin{array}{ll}
\sigma_i\big([\widetilde{H}_{\beta\beta}\ \ \widetilde{H}_{\beta c}]\big)=\sigma_j\big([\widetilde{H}_{\beta\beta}\ \ \widetilde{H}_{\beta c}]\big)\quad \forall\, i,j\in\{|\beta_1|+1,\ldots,|\beta_1|+|\beta_2|\} & \mbox{if $\langle e_{\beta},\overline{u}_{\beta}\rangle= k-k_0$,} \\ [5pt]
\sigma_i\big([\widetilde{H}_{\beta\beta}\ \ \widetilde{H}_{\beta c}]\big) = 0 \quad \forall\, i \in\{|\beta_1|+1,\ldots,|\beta|\} & \mbox{if $\langle e_{\beta},\overline{u}_{\beta}\rangle< k-k_0$,}
\end{array}
\right.
\end{equation}
In fact, it is evident that \eqref{eq:beta2-diag-cri-2} implies that the second equality of \eqref{eq:critical-cone-eq-2} holds. Conversely, consider the following two sub-cases.

{\bf Case 2.1} $\|\overline{S}\|_*=k$, i.e., $\langle e_{\beta},\overline{u}_{\beta}\rangle= k-k_0$. Without loss of generality, assume that $\beta_2\neq\emptyset$, which implies $k-k_0\in\{|\beta_1|+1,\ldots,|\beta_1|+|\beta_2|\}$. Suppose that  there exists $i\in\{|\beta_1|+1,\ldots,|\beta_1|+|\beta_2|\}$ but $i\neq k-k_0$ such that $\sigma_i\big([\widetilde{H}_{\beta\beta}\ \ \widetilde{H}_{\beta c}]\big)>\sigma_{k-k_0}\big([\widetilde{H}_{\beta\beta}\ \ \widetilde{H}_{\beta c}]\big)$ or $\sigma_{k-k_0}\big([\widetilde{H}_{\beta\beta}\ \ \widetilde{H}_{\beta c}]\big)>\sigma_{i}\big([\widetilde{H}_{\beta\beta}\ \ \widetilde{H}_{\beta c}]\big)$. Then, since $0< \overline{u}_{\beta_2}< e_{\beta_2}$ and $\langle e_{\beta},\overline{u}_{\beta} \rangle=k-k_0$, for both cases, we always have 
\begin{eqnarray*}
	&&\sum_{i=1}^{k-k_0}\sigma_i\big([\widetilde{H}_{\beta\beta}\ \ \widetilde{H}_{\beta c}]\big)-\overline{u}_{\beta}^T\sigma\big([\widetilde{H}_{\beta\beta}\ \ \widetilde{H}_{\beta c}]\big) \\ [3pt] &=&\sum_{i=1}^{k-k_0}\sigma_{i}\big([\widetilde{H}_{\beta\beta}\ \ \widetilde{H}_{\beta c}]\big)(1-(\overline{u}_{\beta})_i)-\sum_{i=k-k_0+1}^{|\beta|}(\overline{u}_{\beta})_i\sigma_{i}\big([\widetilde{H}_{\beta\beta}\ \ \widetilde{H}_{\beta c}]\big) \\ [3pt]
	&>& \sum_{i=1}^{k-k_0}\sigma_{k-k_0}\big([\widetilde{H}_{\beta\beta}\ \ \widetilde{H}_{\beta c}]\big)(1-(\overline{u}_{\beta})_i)-\sum_{i=k-k_0+1}^{|\beta|}(\overline{u}_{\beta})_i\sigma_{k-k_0}\big([\widetilde{H}_{\beta\beta}\ \ \widetilde{H}_{\beta c}]\big) \\
	&=& \sigma_{k-k_0}\big([\widetilde{H}_{\beta\beta}\ \ \widetilde{H}_{\beta c}]\big)\left(k-k_0-\sum_{i=1}^{k-k_0}(\overline{u}_{\beta})_i-\sum_{i=k-k_0+1}^{|\beta|}(\overline{u}_{\beta})_i\right)=0,
\end{eqnarray*}
which implies that the second equality of \eqref{eq:critical-cone-eq-2} does not hold, which contradicts the assumption. Therefore, we know that $H\in{\cal C}(X;\partial\,\theta(\overline{X}))$ if and only if (ii) holds. 

{\bf Case 2.2} $\|\overline{S}\|_*<k$, i.e., $\langle e_{\beta},\overline{u}_{\beta} \rangle< k-k_0$. We know that $\beta_2\cup\beta_3\neq\emptyset$ and $k-k_0\in\{|\beta_1|+1,\ldots,|\beta|\}$. Suppose that \eqref{eq:beta2-diag-cri-2} does not hold. Then, we know that either there exists $i\in\{|\beta_1|+1,\ldots,|\beta|\}$ such that $i<k-k_0$ and $\sigma_{i}\big([\widetilde{H}_{\beta\beta}\ \ \widetilde{H}_{\beta c}]\big)>\sigma_{k-k_0}\big([\widetilde{H}_{\beta\beta}\ \ \widetilde{H}_{\beta c}]\big)=0$ or $\sigma_{k-k_0}\big([\widetilde{H}_{\beta\beta}\ \ \widetilde{H}_{\beta c}]\big)>0$. For the case that $\sigma_{i}\big([\widetilde{H}_{\beta\beta}\ \ \widetilde{H}_{\beta c}]\big)>\sigma_{k-k_0}\big([\widetilde{H}_{\beta\beta}\ \ \widetilde{H}_{\beta c}]\big)=0$, since $0< \overline{u}_{\beta_2}< e_{\beta_2}$, we have 
\begin{eqnarray*}
	&&\sum_{i=1}^{k-k_0}\sigma_i\big([\widetilde{H}_{\beta\beta}\ \ \widetilde{H}_{\beta c}]\big)-\overline{u}_{\beta}^T\sigma\big([\widetilde{H}_{\beta\beta}\ \ \widetilde{H}_{\beta c}]\big) \\ [3pt] &=&\sum_{i=1}^{k-k_0}\sigma_{i}\big([\widetilde{H}_{\beta\beta}\ \ \widetilde{H}_{\beta c}]\big)(1-(\overline{u}_{\beta})_i)-\sum_{i=k-k_0+1}^{|\beta|}(\overline{u}_{\beta})_i\sigma_{i}\big([\widetilde{H}_{\beta\beta}\ \ \widetilde{H}_{\beta c}]\big) \\ [3pt]
	&>& \sum_{i=1}^{k-k_0}\sigma_{k-k_0}\big([\widetilde{H}_{\beta\beta}\ \ \widetilde{H}_{\beta c}]\big)(1-(\overline{u}_{\beta})_i)-\sum_{i=k-k_0+1}^{|\beta|}(\overline{u}_{\beta})_i\sigma_{k-k_0}\big([\widetilde{H}_{\beta\beta}\ \ \widetilde{H}_{\beta c}]\big)=0.
\end{eqnarray*}
For the case that $\sigma_{k-k_0}\big([\widetilde{H}_{\beta\beta}\ \ \widetilde{H}_{\beta c}]\big)>0$, since $\langle e_{\beta},\overline{u}_{\beta}\rangle< k-k_0$, we obtain that 
\begin{eqnarray*}
	&&\sum_{i=1}^{k-k_0}\sigma_i\big([\widetilde{H}_{\beta\beta}\ \ \widetilde{H}_{\beta c}]\big)-\overline{u}_{\beta}^T\sigma\big([\widetilde{H}_{\beta\beta}\ \ \widetilde{H}_{\beta c}]\big) \\ [3pt] &=&\sum_{i=1}^{k-k_0}\sigma_{i}\big([\widetilde{H}_{\beta\beta}\ \ \widetilde{H}_{\beta c}]\big)(1-(\overline{u}_{\beta})_i)-\sum_{i=k-k_0+1}^{|\beta|}(\overline{u}_{\beta})_i\sigma_{i}\big([\widetilde{H}_{\beta\beta}\ \ \widetilde{H}_{\beta c}]\big) \\ [3pt]
	&\ge& \sum_{i=1}^{k-k_0}\sigma_{k-k_0}\big([\widetilde{H}_{\beta\beta}\ \ \widetilde{H}_{\beta c}]\big)(1-(\overline{u}_{\beta})_i)-\sum_{i=k-k_0+1}^{|\beta|}(\overline{u}_{\beta})_i\sigma_{k-k_0}\big([\widetilde{H}_{\beta\beta}\ \ \widetilde{H}_{\beta c}]\big) \\
	&=& \sigma_{k-k_0}\big([\widetilde{H}_{\beta\beta}\ \ \widetilde{H}_{\beta c}]\big)\left(k-k_0-\sum_{i=1}^{k-k_0}(\overline{u}_{\beta})_i-\sum_{i=k-k_0+1}^{|\beta|}(\overline{u}_{\beta})_i\right)>0.
\end{eqnarray*}
Therefore, for both cases, we always conclude that the second equality in \eqref{eq:critical-cone-eq-2} does not hold, which contradicts the assumption. Therefore, we know that $H\in{\cal C}(X;\partial\,\theta(\overline{X}))$ if and only if (iii) holds.  $\hfill\Box$

\vskip 10 true pt

For the given $\overline{S}\in\partial\,\theta(\overline{X})$, let ${\rm aff}({\cal C}(X;\partial\,\theta(\overline{X}))$ be the affine hull of the critical cone ${\cal C}(X;\partial\,\theta(\overline{X}))$, i.e., the smallest affine space containing ${\cal C}(X;\partial\,\theta(\overline{X}))$. Note that it follows from \eqref{eq:eqiv-critical-2} that $0\in{\cal C}(X;\partial\,\theta(\overline{X})$. It is easy to see (cf. e.g., \cite[Theorem 2.7]{Rockafellar70}) that ${\rm aff}({\cal C}(X;\partial\,\theta(\overline{X}))={\cal C}(X;\partial\,\theta(\overline{X}))-{\cal C}(X;\partial\,\theta(\overline{X}))$.  Therefore, by Proposition \ref{prop:critical-cone}, one can easily derive the following proposition on the characterization of ${\rm aff}({\cal C}(X;\partial\,\theta(\overline{X}))$. For simplicity, we omit the detail proof here.
\begin{proposition}
	 Suppose that $(\overline{X},\overline{S})\in\Re^{m\times n}\times \Re^{m\times n}$ is a solution of the GE \eqref{eq:GE}. Let $X=\overline{X}+\overline{S}$ admit the SVD \eqref{eq:SVD-X}. Then, $H\in {\rm aff}({\cal C}(X;\partial\,\theta(\overline{X}))$ if and only if $H$ satisfies the following conditions.
	 \begin{itemize}
	 	\item[(i)] If $\sigma_k(\overline{X})>0$, then there exists some $\tau\in\Re$ such that
	 	\begin{equation*}\label{eq:affC-Case1}
	 	S(\overline{U}_{\beta}^TH\overline{V}_{\beta})= \left[
	 	\begin{array}{ccc}
	 	S(\overline{U}_{\beta_1}^TH\overline{V}_{\beta_1}) & 0 & 0\\
	 	0 & \tau I_{|\beta_2|} & 0  \\
	 	0 & 0 & S(\overline{U}_{\beta_3}^TH\overline{V}_{\beta_3}) 
	 	\end{array}
	 	\right].
	 	\end{equation*}
	 	\item[(ii)] If $\sigma_k(\overline{X})=0$ and $\|\overline{S}\|_*=k$, then there exists some $\tau\in\Re$ such that
	 	\begin{equation*}\label{eq:affC-Case2.1}
	 	\left[\overline{U}_{\beta}^{T}H\overline{V}_{\beta}\ \ \overline{U}_{\beta}^{T}H\overline{V}_{2}\right]=\left[\begin{array}{cccc}
	 	S(\overline{U}_{\beta_1}^{T}H\overline{V}_{\beta_1}) & 0 & 0 & 0 \\
	 	0 & \tau I_{|\beta_2|} &  0 & 0 \\
	 	0  & 0 &  \overline{U}_{b}^{T}H\overline{V}_{b} & \overline{U}_{b}^{T}H\overline{V}_{2} \end{array}\right].
	 	\end{equation*}
	 	\item[(iii)] If $\sigma_k(\overline{X})=0$ and $\|\overline{S}\|_*<k$, then 
	 	\begin{equation*}\label{eq:affC-Case2.2}
	 	\left[\overline{U}_{\beta}^{T}H\overline{V}_{\beta}\ \ \overline{U}_{\beta}^{T}H\overline{V}_{2}\right]=\left[\begin{array}{cccc}
	 	S(\overline{U}_{\beta_1}^{T}H\overline{V}_{\beta_1}) & 0 & 0 & 0 \\
	 	0 & 0 &  0 & 0 \\
	 	0  & 0 &  0 & 0 \end{array}\right].
	 	\end{equation*}
	 \end{itemize}
\end{proposition}

%Moreover, by comparing with the characterizations of Clarke's generalized Jacobian of the proximal mapping  ${\rm Pr}_{\theta}$ (Proposition \ref{prop:Bsubdiff-Proj-knorm}), we obtain the following simple observation.
%
%
%
%\begin{lemma}\label{lemma:V-in-affC}
%	Let $X=\overline{X}+\overline{S}$ with $\overline{S}\in\partial\,\theta(\overline{X})$.  For any ${\cal V} \in\partial\,{\rm Pr}_{\theta}(X)$, we have
%	\[
%	 {\cal V}(H)\in {\rm aff}\left({\cal C}(X;\partial\,\theta(\overline{X}))\right),\quad H\in\Re^{m\times n}.
%	\]
%\end{lemma}
%\noindent{\bf Proof.} See the Appendix. $\hfill \Box$
%
%\vskip 10 true pt

Next, consider the dual GE \eqref{eq:GE-dual}. The critical cone of ${\cal K}^{\circ}$ at $(X,t)=(\overline{X}+\overline{S},\theta(\overline{X})-1)\in\Re^{m\times n}\times\Re$, associated with the complementarity problem in \eqref{eq:GE=CP}, is defined as
\begin{equation}\label{eq:def-critical cone-Ko}
{\cal C}\left((X,t);{\cal K}^{\circ}\right)={\cal T}_{{\cal K}^{\circ}}(\overline{S},-1)\cap(\overline{X},\theta(\overline{X}))^{\perp}.
\end{equation}
Thus, we know from \eqref{eq:T-epi-f} that 
\begin{equation}\label{eq:eqiv-critical-1-dual}
(\widehat{H},\tau)\in{\cal C}\left((X,t);{\cal K}^{\circ}\right) \Longleftrightarrow 
\left\{
\begin{array}{l}
\widehat{H}=H-\overline{U}\left[\begin{array}{cc}
\tau I_k & 0 \\
0 & 0
\end{array}\right]\overline{V}^T, \\ [3pt]
H\in{\cal C}(X;\partial\,\theta^*(\overline{S})), 
\end{array}
\right. 
\end{equation}
where ${\cal C}(X;\partial\,\theta^*(\overline{S}))\subseteq\Re^{m\times n}$ is defined by
\begin{equation}\label{eq:def-critical-cone-dual}
{\cal C}(X;\partial\,\theta^*(\overline{S})):=\left\{
\begin{array}{ll}
\big\{H\in\Re^{m\times n}\mid \vartheta'(\overline{S};H)\le \langle \overline{X},H\rangle=0 \big\} & \mbox{if $\vartheta(\overline{S})=1$,} \\ [3pt]
\Re^{m\times n} & \mbox{if $\vartheta(\overline{S})<1$.}
\end{array} \right. 
\end{equation}
We call ${\cal C}(X;\partial\,\theta^*(\overline{S}))$ the critical cone of $\partial\,\theta^*(\overline{S})={\cal N}_{{\cal B}_{(k)^*}}(\overline{S})$ at $X=\overline{X}+\overline{S}$, associated with the dual GE in \eqref{eq:GE-dual}. The following characterization of the critical cone ${\cal C}(X;\partial\,\theta^*(\overline{S}))$ can be obtain similarly as that of ${\cal C}(X;\partial\,\theta(\overline{X}))$. For simplicity, we omit the detail proof here.

\begin{proposition}\label{prop:critical-cone-dual}
	Suppose that $(\overline{X},\overline{S})\in\Re^{m\times n}\times \Re^{m\times n}$ is a solution of the dual GE \eqref{eq:GE-dual}. Assume that  $\vartheta(\overline{S})=1$. Let $X=\overline{X}+\overline{S}$ admit the SVD \eqref{eq:SVD-X}.  Then, 
	\begin{equation*}
	H\in{\cal C}(X;\partial\,\theta^*(\overline{S}))\quad \Longleftrightarrow \quad \vartheta'(\overline{S};H)=\langle \overline{X},H\rangle=0,
	\end{equation*}
	which is equivalent to the following conditions. 
	\begin{itemize}
		\item[(i)] If $\sigma_k(\overline{X})>0$, then ${\rm tr}(\overline{U}_{\beta}^TH\overline{V}_{\beta})=0$,
		\[
		S(\overline{U}_{\alpha\cup\beta_1}^TH\overline{V}_{\alpha\cup\beta_1})= \left[
		\begin{array}{cc}
		0 & 0\\
		0 & S(\overline{U}_{\beta_1}^TH\overline{V}_{\beta_1}) 
		\end{array}
		\right] \quad {\rm with} \quad S(\overline{U}_{\beta_1}^TH\overline{V}_{\beta_1})\preceq 0 
		\]
		and
		\[
		\left[\overline{U}_{\beta_3\cup\gamma}^{T}H\overline{V}_{\beta_3\cup\gamma}\ \ \overline{U}_{\beta_3\cup\gamma}^{T}H\overline{V}_{2}\right]=\left[
		\begin{array}{ccc}
		S(\overline{U}_{\beta_3}^TH\overline{V}_{\beta_3}) & 0 & 0\\
		0 & 0 & 0 
		\end{array}
		\right] \quad {\rm with} \quad S(\overline{U}_{\beta_3}^TH\overline{V}_{\beta_3})\succeq 0.
		\]
		\item[(ii)] If $\sigma_k(\overline{X})=0$ and $\|\overline{S}\|_*<k$, then
		\[
		S(\overline{U}_{\alpha\cup\beta_1}^TH\overline{V}_{\alpha\cup\beta_1})= \left[
		\begin{array}{cc}
		0 & 0\\
		0 & S(\overline{U}_{\beta_1}^TH\overline{V}_{\beta_1}) 
		\end{array}
		\right] \quad {\rm with} \quad S(\overline{U}_{\beta_1}^TH\overline{V}_{\beta_1})\preceq 0.
		\] 
		\item[(iii)] If $\sigma_k(\overline{X})=0$ and $\|\overline{S}\|_*=k$, then ${\rm tr}(\overline{U}_{\beta_1\cup\beta_2}^TH\overline{V}_{\beta_1\cup\beta_2})+\left\|\big[\overline{U}_{b}^{T}H\overline{V}_{b}\ \ \overline{U}_{b}^{T}H\overline{V}_{2}\big]\right\|_*\le 0$,
		\[
		S(\overline{U}_{\alpha\cup\beta_1}^TH\overline{V}_{\alpha\cup\beta_1})= \left[
		\begin{array}{cc}
		0 & 0\\
		0 & S(\overline{U}_{\beta_1}^TH\overline{V}_{\beta_1}) 
		\end{array}
		\right] \quad {\rm with} \quad S(\overline{U}_{\beta_1}^TH\overline{V}_{\beta_1})\preceq 0.
		\] 
	\end{itemize}
	
\end{proposition}

For the given $\overline{X}\in\partial\,\theta^*(\overline{S})$, let ${\rm aff}({\cal C}(X;\partial\,\theta^*(\overline{S}))$ be the affine hull of the critical cone ${\cal C}(X;\partial\,\theta^*(\overline{S}))$. Therefore, by Proposition \ref{prop:critical-cone-dual}, we obtain the following characterization of ${\rm aff}({\cal C}(X;\partial\,\theta^*(\overline{S}))$.
\begin{proposition}
	Suppose that $(\overline{X},\overline{S})\in\Re^{m\times n}\times \Re^{m\times n}$ is a solution of the dual GE \eqref{eq:GE-dual}. Assume that  $\vartheta(\overline{S})=1$. Let $X=\overline{X}+\overline{S}$ admit the SVD \eqref{eq:SVD-X}. Then, $H\in {\rm aff}({\cal C}(X;\partial\,\theta^*(\overline{S}))$ if and only if $H$ satisfies the following conditions.
	\begin{itemize}
		\item[(i)] If $\overline{\sigma}_k>0$, then 
		\begin{equation}\label{eq:affC*-Case1.1}
		S(\overline{U}_{\alpha\cup\beta_1}^TH\overline{V}_{\alpha\cup\beta_1})= \left[
		\begin{array}{cc}
		0 & 0\\
		0 & S(\overline{U}_{\beta_1}^TH\overline{V}_{\beta_1}) 
		\end{array}
		\right], \quad {\rm tr}(\overline{U}_{\beta}^TH\overline{V}_{\beta})=0.
		\end{equation}
		and
		\begin{equation}\label{eq:affC*-Case1.2}
		\left[\overline{U}_{\beta_3\cup\gamma}^{T}H\overline{V}_{\beta_3\cup\gamma}\ \ \overline{U}_{\beta_3\cup\gamma}^{T}H\overline{V}_{2}\right]=\left[
		\begin{array}{ccc}
		S(\overline{U}_{\beta_3}^TH\overline{V}_{\beta_3}) & 0 & 0\\
		0 & 0 & 0 
		\end{array}
		\right].
		\end{equation}
		\item[(ii)] If $\overline{\sigma}_k=0$, then  
		\begin{equation}\label{eq:affC*-Case2}
		S(\overline{U}_{\alpha\cup\beta_1}^TH\overline{V}_{\alpha\cup\beta_1})= \left[
		\begin{array}{cc}
		0 & 0\\
		0 & S(\overline{U}_{\beta_1}^TH\overline{V}_{\beta_1}) 
		\end{array}
		\right].
		\end{equation}
	\end{itemize}
\end{proposition}

\section{The second order analysis}\label{sect:second order tangent sets}
In this section, we shall study another important variational property of the Ky Fan $k$-norm $\theta=\|\cdot\|_{(k)}$, i.e., the conjugate function of the parabolic second order directional derivative of $\theta$, which equals to the support function of the second order tangent set of the epigraph of $\theta$. This conjugate function is closely related to the second order optimality conditions of the problem \eqref{eq:P-g}. 

For the given $(\overline{X},\theta(\overline{X}))\in{\cal K}$, let ${\cal T}^{i,2}_{\cal K}\left((\overline{X},\theta(\overline{X}));(H,\tau)\right)$ and ${\cal T}^{2}_{\cal
K}\left((\overline{X},\theta(\overline{X}));(H,\tau)\right)$ be the inner and outer second order tangent sets \cite[Definition 3.28]{BShapiro00} to  ${\cal K}$ at $(\overline{X},\theta(\overline{X}))\in{\cal K}$ along the direction $(H,\tau)\in{\cal T}_{{\cal K}}(\overline{X},\theta(\overline{X}))$, respectively, i.e.,
\[
{\cal T}^{i,2}_{\cal K}\left((\overline{X},\theta(\overline{X}));(H,\tau)\right):=\liminf_{\rho\downarrow 0}\frac{{\cal K}-(\overline{X},\theta(\overline{X}))-\rho(H,\tau)}{\frac{1}{2}\rho^{2}}
\]
and
\[
{\cal T}^{2}_{\cal
	K}\left((\overline{X},\theta(\overline{X}));(H,\tau)\right):=\limsup_{\rho\downarrow 0}\frac{{\cal K}-(\overline{X},\theta(\overline{X}))-\rho(H,\tau)}{\frac{1}{2}\rho^{2}},
\]
where ``$\limsup$'' and ``$\liminf$'' are the Painlev\'{e}-Kuratowski outer and inner limit for sets (cf. \cite[Definition 4.1]{RWets98}). 

For ${\cal T}^2_{\cal K}:={\cal T}^{i,2}_{\cal K}\left((\overline{X},\theta(\overline{X}));(H,\tau)\right)$ or ${\cal T}^{2}_{\cal K}\left((\overline{X},\theta(\overline{X}));(H,\tau)\right)$, since ${\cal K}$ is convex, we know from \cite[Proposition 3.34, (3.62) \& (3.63)]{BShapiro00}  that for any $(\overline{X},\theta(\overline{X}))\in{\cal K}$ and $(H,\tau)\in{\cal T}_{\cal K}(\overline{X},\theta(\overline{X}))$, 
\begin{equation}\label{eq:second-order-tangent-inclusions}
{\cal T}^2_{\cal K}+{\cal T}_{{\cal T}_{\cal K}(\overline{X},\theta(\overline{X}))}(H,\tau)\subseteq {\cal T}^2_{\cal K}\subseteq{\cal T}_{{\cal T}_{\cal K}(\overline{X},\theta(\overline{X}))}(H,\tau), 
\end{equation}
where ${\cal T}_{{\cal T}_{\cal K}(\overline{X},\theta(\overline{X}))}(H,\tau)$ is the tangent cone of ${\cal T}_{\cal K}(\overline{X},\theta(\overline{X}))$ at $(H,\tau)$.  
%\begin{proposition}\label{prop:second-order-tangent-inclusions}
%Let $C$ be a convex set. Then, for any $x\in C$, $h\in{\cal T}_{C}(x)$, the following inclusions hold:
%\[
%{\cal T}^{i,2}_{C}(x,h)+{\cal T}_{{\cal T}_{C}(x)}(h)\subseteq{\cal T}^{i,2}_{C}(x,h)\subseteq{\cal T}_{{\cal T}_{C}(x)}(h),
%\]
%\[
%{\cal T}^{2}_{C}(x,h)+{\cal T}_{{\cal T}_{C}(x)}(h)\subseteq{\cal T}^{2}_{C}(x,h)\subseteq{\cal T}_{{\cal T}_{C}(x)}(h),
%\]
%where ${\cal T}_{{\cal T}_{C}(x)}(h)$ is the tangent cone of ${\cal T}_{C}(x)$ at $h\in{\cal T}_{C}(x)$.
%\end{proposition}
For any given $(H,\tau)\in{\cal T}_{{\cal K}}(\overline{X},\theta(\overline{X}))$, let us consider the following two cases.

{\bf Case 1.} $\sum_{i=1}^{k}\sigma'_{i}(\overline{X};H)=\tau$, i.e., $(H,\tau)\in{\rm bd}\,{\cal T}_{\cal K}(\overline{X},\theta(\overline{X}))$. 
Since ${\rm int}\,{\cal K}\neq \emptyset$ and the continuous convex function $\theta=\|\cdot\|_{(k)}$ is
(parabolically) second order directionally differentiable, we know from
\cite[Proposition 3.30]{BShapiro00} that
\[
{\cal T}^{i,2}_{\cal K}\left((\overline{X},\theta(\overline{X}));(H,\tau)\right)={\cal T}^{2}_{\cal K}\left((\overline{X},\theta(\overline{X}));(H,\tau)\right)={\rm epi}\, \theta''(\overline{X};H,\cdot),
\]
where ${\rm epi}\,\theta''(\overline{X};H,\cdot)$ is the epigraph of the (parabolic) second order directional derivative of $\theta$ at $\overline{X}$ along the direction $H$, which is convex and given by
\begin{equation}\label{eq:bd2-secondordertangent}
{\rm epi}\,\theta''(\overline{X};H,\cdot):=\Big\{ (W,\eta)\in\Re^{m\times n}\times\Re\,|\,
\sum_{i=1}^{k}\sigma''_{i}(\overline{X};H,W)\leq \eta \Big\}.
\end{equation}
 
{\bf Case 2.} $\sum_{i=1}^{k}\sigma'_{i}(\overline{X};H)<\tau$,
i.e., $(H,\tau)\in{\rm int}\,{\cal T}_{\cal K}(\overline{X},\theta(\overline{X}))$. Since ${\cal T}_{{\cal T}_{\cal K}(\overline{X},\theta(\overline{X}))}(H,\tau)=\Re\times\Re^{m\times n}$, we know from \eqref{eq:second-order-tangent-inclusions} that
\begin{equation}\label{eq:bd1-secondordertangent}
{\cal T}^{i,2}_{\cal K}\left((\overline{X},\theta(\overline{X}));(H,\tau)\right)={\cal T}^{2}_{\cal K}\left((\overline{X},\theta(\overline{X}));(H,\tau)\right)=\Re\times\Re^{m\times n}.
\end{equation}
Therefore, we may denote ${\cal T}^{2}_{\cal K}\left((\overline{X},\theta(\overline{X}));(H,\tau)\right)$ the second order tangent set to ${\cal K}$ at $(\overline{X},\theta(\overline{X}))$ along the direction $(H,\tau)\in{\cal T}_{{\cal K}}(\overline{X},\theta(\overline{X}))$.

\vskip 10 true pt

Next, we shall provide the explicit formula of the support function of the second order tangent set ${\cal T}^{2}_{\cal K}\left((\overline{X},\theta(\overline{X}));(H,\tau)\right)$.  Let $(\overline{X},\theta(\overline{X}))\in {\cal K}$ be fixed. For any $(H,\tau)\in{\cal T}_{{\cal K}}(\overline{X},\theta(\overline{X}))$, denote ${\cal T}^2(H,\tau):={\cal T}^{2}_{\cal K}\left((\overline{X},\theta(\overline{X}));(H,\tau)\right)$.  Consider the support function $\delta^{*}_{{\cal T}^2(\tau,H)}(\cdot,\cdot):\Re\times\Re^{m\times n}\to(-\infty,\infty]$, i.e.,
\[
\delta^{*}_{{\cal T}^2(H,\tau)}(S,\zeta)=\sup\left\{\langle S, W\rangle+\zeta\eta \mid
(W,\eta)\in{\cal T}^2(H,\tau)\right\},\quad (S,\zeta)\in\Re^{m\times n}\times\Re.
\]

\begin{claim}\label{claim1}
 $\delta^{*}_{{\cal T}^2(H,\tau)}(S,\zeta)\equiv \infty$ if  $(S,\zeta)\notin \big({\cal T}_{{\cal T}_{\cal K}(\overline{X},\theta(\overline{X}))}(H,\tau)\big)^{\circ}$.
\end{claim} 
\noindent{\bf Proof.} Let $(S,\zeta)\notin \big({\cal T}_{{\cal T}_{\cal K}(\overline{X},\theta(\overline{X}))}(H,\tau)\big)^{\circ}$ be arbitrarily given. Since ${\cal T}_{{\cal T}_{\cal K}(\overline{X},\theta(\overline{X}))}(H,\tau)$ is nonempty, we may assume that
there exists $(W^{\circ},\eta^{\circ})\in {\cal T}_{{\cal T}_{\cal K}(\overline{X},\theta(\overline{X}))}(H,\tau)$ such that
\[
\langle (S,\zeta), (W^{\circ},\eta^{\circ}) \rangle >0.
\] 
Fix any $(\widetilde{\eta},\widetilde{W})\in{\cal T}^2(H,\tau)$.
By \eqref{eq:second-order-tangent-inclusions}, we have for any $\rho>0$,  
\[
\rho(W^{\circ},\eta^{\circ})+(\widetilde{W},\widetilde{\eta})\in {\cal T}_{{\cal T}_{\cal K}(\overline{X},\theta(\overline{X}))}(H,\tau)+{\cal T}^2(H,\tau)\subseteq {\cal T}^2(H,\tau).
\] 
Therefore, we know that
\[
\rho\langle (S,\zeta), (W^{\circ},\eta^{\circ}) \rangle+\langle (S,\zeta),
(\widetilde{W},\widetilde{\eta}) \rangle \leq \delta^{*}((S,\zeta)\mid {\cal T}^2(H,\tau)).
\] 
Since $\langle (S,\zeta), (W^{\circ},\eta^{\circ}) \rangle >0$ and $\rho>0$ can be arbitrarily large, we conclude that $\delta^{*}_{{\cal T}^2(H,\tau)}(S,\zeta)\equiv\infty$ for any $(S,\zeta)\notin \big({\cal T}_{{\cal T}_{\cal K}(\overline{X},\theta(\overline{X}))}(H,\tau)\big)^{\circ}$. $\hfill \Box$

\vskip 10 true pt

Since ${\cal K}$ is a closed convex cone in $\Re^{m\times n}\times\Re$, it can be verified easily that
\[
{\cal K}\subseteq {\cal T}_{{\cal K}}(\overline{X},\theta(\overline{X}))\subseteq {\cal T}_{{\cal T}_{\cal K}(\overline{X},\theta(\overline{X}))}(H,\tau).
\] In particular, we have $\pm (\overline{X},\theta(\overline{X}))\in {\cal T}_{{\cal K}}(\overline{X},\theta(\overline{X})) \subseteq {\cal T}_{{\cal T}_{\cal K}(\overline{X},\theta(\overline{X}))}(H,\tau)$ and $\pm (H,\tau)\in {\cal T}_{{\cal T}_{\cal K}(\overline{X},\theta(\overline{X}))}(H,\tau)$. Therefore, we know from the definition of the polar cone that if $(S,\zeta)\in\big({\cal T}_{{\cal T}_{\cal K}(\overline{X},\theta(\overline{X}))}(H,\tau)\big)^{\circ}$, then
\begin{equation}\label{eq:multipler-condition}
(S,\zeta)\in{\cal K}^{\circ},\quad \left\langle (S,\zeta), (\overline{X},\theta(\overline{X}))\right\rangle = 0 \quad
{\rm and}  \quad \left\langle (S,\zeta), (H,\tau)\right\rangle = 0.
\end{equation} 
Hence, by Claim \ref{claim1}, we only need to consider the point $(S,\zeta)\in\Re\times\Re^{m\times n}$ satisfying the
condition (\ref{eq:multipler-condition}), since otherwise $\delta^{*}_{{\cal T}^{2}}(S,\zeta)\equiv \infty$. Moreover, instead of considering the general  $S\in\Re^{m\times n}$, we only consider the point $\overline{S}$ such that $(\overline{X},\overline{S})\in\Re^{m\times n}\times \Re^{m\times n}$ satisfying the GE \eqref{eq:GE}, i.e., $\overline{S}\in\partial\,\theta(\overline{X})$, which is equivalent to the complementarity problem in \eqref{eq:GE=CP}. 

On the other hand, by the definition of the critical cone \eqref{eq:def-critical cone-K} of ${\cal K}$, it is evident that the given point $(\overline{S},-1)$ satisfies the condition (\ref{eq:multipler-condition}) if and only if $(H,\tau)\in{\cal C}((X,t);{\cal K})$ with $(X,t)=(\overline{X},\theta(\overline{X}))+(\overline{S},-1)$. Thus, by \eqref{eq:eqiv-critical-1} and \eqref{eq:eqiv-critical-2}, we know that $(\overline{S},-1)$ satisfies the condition (\ref{eq:multipler-condition}) if and only if $H\in{\cal C}(X;\partial\,\theta(\overline{X}))$ (defined by \eqref{eq:def-critical-cone}) and $\tau=\langle \overline{S},H\rangle=\sum_{i=1}^{k}\sigma'_{i}(\overline{X};H)$. Hence, we know from \eqref{eq:bd2-secondordertangent} that
\begin{equation}\label{eq:second-order-tangent-used}
{\cal T}^2(H):={\cal T}^2(H,\tau)=\Big\{ (W,\eta)\in\Re^{m\times n}\times\Re\,|\,
\sum_{i=1}^{k}\sigma''_{i}(\overline{X};H,W)\leq \eta \Big\},
\end{equation}
where for each $i$, the  second order directional derivative  $\sigma''_i(\overline{X},H,W)$ is given by Proposition \ref{prop:second-directional-diff-singlevalue}.  Let $X=\overline{X}+\overline{S}$ admit the SVD \eqref{eq:SVD-X}. Let $a_{1},\ldots,a_{r}$ be the index sets defined by (\ref{eq:ak-nonsymmetric}) with respect to $X$. Denote $\overline{\sigma}=\sigma(\overline{X})$ and $\overline{u}=\sigma(\overline{S})$.  Consider the following two cases.  

{\bf Case 1.} $\overline{\sigma}_k>0$. Let $\alpha$, $\beta$ and $\gamma$ be the index sets defined by \eqref{eq:def-alpha-beta-gamma-case1} and $\beta_1$, $\beta_2$ and $\beta_3$ be the index sets defined by \eqref{eq:def-beta123-case1}.  Let $r_{0} \le r_{1}\in\{0,1,\ldots,r+1\}$, $r_{0} \le\widetilde{r}_0\le r_{0}+1$ and $r_{1}-1 \le\widetilde{r}_1\le r_1$ be the integers  such that \eqref{eq:widetilde_r-case1} holds. For each $l\in\{1,\ldots,r_0\}$, since $\overline{\sigma}_i=\overline{\sigma}_{i'}$ for any $i,i'\in a_l$, we use  $\overline{\nu}_l$ to denote the common value. By \eqref{eq:eqiv-critical-1} and \eqref{eq:eqiv-critical-2} , we know that there exists an orthogonal matrix $R\in{\cal O}^{|\beta|}$ such that \eqref{eq:simultaneous-ordered} holds, i.e., ${\rm Diag}(\overline{u}_{\beta})$ and $S(\overline{U}^T_{\beta}H\overline{V}_{\beta})$ admit a simultaneous ordered eigenvalue decomposition. Therefore, $R$ has the block diagonal structure \eqref{eq:R-block-diag-1}.  Hence, we know from the part (i) of Proposition \ref{prop:second-directional-diff-singlevalue} that $(W,\eta)\in{\cal T}^2(H)$ if and only if 
\begin{eqnarray}
\sum_{i=1}^{k}\sigma''_i(\overline{X};H,W) &=&\sum_{l=1}^{r_0}{\rm tr}\,(S(\overline{U}^T_{a_{l}}W\overline{V}_{a_{l}}))-2\sum_{l=1}^{r_0}{\rm tr}\left(\Omega_{a_l}(\overline{X},H)\right)\nonumber \\ 
&& +{\rm tr}\left(R_1^{T}\big(S(\overline{U}^T_{\beta}W\overline{V}_{\beta})-2\Omega_{\beta}(\overline{X},H)\big) R_1\right)\nonumber \\ 
&&+\sum_{i=|\beta_1|+1}^{k-k_0}\lambda_i\left(R_2^T\big(S(\overline{U}^T_{\beta}W\overline{V}_{\beta})-2\Omega_{\beta}(\overline{X},H)\big)R_2\right) \leq \eta,\label{eq:eta-WinT2}
\end{eqnarray}
where $\Omega_{a_l}(\overline{X},H)\in{\cal S}^{m}$, $l=1,\ldots,r_0$ and $\Omega_{\beta}(\overline{X},H)\in{\cal S}^{m}$ are given by \eqref{eq:def-Omega-dd} with respect to $\overline{X}$, $R_1\in{\cal O}(S(\overline{U}_{\beta_1}^TH\overline{V}_{\beta_1}))$ and $R_2\in{\cal O}(S(\overline{U}_{\beta_2}^TH\overline{V}_{\beta_2}))$.  
Meanwhile, since $\overline{u}_{\alpha}=e_{\alpha}$, we have for any $(W,\eta)\in{\cal T}^2(H)$, 
\begin{eqnarray}
	&&-\eta+\langle \overline{S}, W \rangle = -\eta+ \left\langle \overline{U}^T\overline{S}\overline{V}, \overline{U}^TW\overline{V} \right\rangle \nonumber\\ [3pt]
	&=& -\eta +\left\langle \left[\begin{array}{cc}
		{\rm Diag}(\overline{u}_{\alpha}) & 0 \\
		0 & {\rm Diag}(\overline{u}_\beta)
	\end{array}\right], \left[\begin{array}{cc}
	S(\overline{U}_{\alpha}^TW\overline{V}_{\alpha}) & 0 \\
	0 & S(\overline{U}_{\beta}^TW\overline{V}_{\beta})
\end{array}\right] \right\rangle \nonumber\\ [3pt]
	&=& -\eta+\sum_{l=1}^{r_0}{\rm tr}\,(S(\overline{U}^T_{a_{l}}W\overline{V}^T_{a_{l}}))+\left\langle {\rm Diag}(\overline{u}_\beta)
	,  S(\overline{U}^T_{\beta}W\overline{V}^T_{\beta})  \right\rangle\nonumber\\ [3pt]
	&=& \Xi(W,\eta)  +2\sum_{l=1}^{r_0}{\rm tr}\left(\Omega_{a_l}(\overline{X},H)\right)+\left\langle {\rm Diag}(\overline{u}_\beta), 2\Omega_{\beta}(\overline{X},H)\right\rangle ,\label{eq:support-product}
\end{eqnarray}
where
\begin{eqnarray}
\Xi(W,\eta)&=&-\eta+ \sum_{l=1}^{r_0}{\rm tr}\,(S(\overline{U}^T_{a_{l}}W\overline{V}^T_{a_{l}}))-2\sum_{l=1}^{r_0}{\rm tr}\left(\Omega_{a_l}(\overline{X},H)\right) \nonumber\\ 
&&+\left\langle {\rm Diag}(\overline{u}_\beta),  S(\overline{U}^T_{\beta}W\overline{V}^T_{\beta})-2\Omega_{\beta}(\overline{X},H)\right\rangle.\label{eq:Delta=case1}
\end{eqnarray}

Next, we shall show that
\begin{equation}\label{eq:Xi=0-1}
\max \left\{\Xi(W,\eta) \,|\, (W,\eta)\in{\cal T}^2(H)\right\}=0.
\end{equation}
In fact, since $0\le\overline{u}_{\beta}\le e_\beta$ and $\langle e_{\beta},\overline{u}\rangle=k-k_0$, we know from Lemma \ref{lem:Fan} (Fan's inequality) that the last term of \eqref{eq:Delta=case1} satisfies 
\begin{eqnarray*}
&& \left\langle {\rm Diag}(\overline{u}_\beta),  S(\overline{U}^T_{\beta}W\overline{V}^T_{\beta})-2\Omega_{\beta}(\overline{X},H)\right\rangle=\left\langle {\rm Diag}(\overline{u}_\beta), R^{T}\left(S(\overline{U}^T_{\beta}W\overline{V}^T_{\beta})-2\Omega_{\beta}(\overline{X},H)\right)R \right\rangle\nonumber\\ 
&\leq&{\rm tr}\left(R_1^{T}\big(S(\overline{U}^T_{\beta}W\overline{V}_{\beta})-2\Omega_{\beta}(\overline{X},H)\big)R_1\right) +  \left\langle \overline{u}_{\beta_2}, \lambda\left(R_2^T\big(S(\overline{U}^T_{\beta}W\overline{V}_{\beta})-2\Omega_{\beta}(\overline{X},H)\big)R_2\right)\right\rangle\nonumber\\ 
&\leq&{\rm tr}\left(R_1^{T}\big(S(\overline{U}^T_{\beta}W\overline{V}_{\beta})-2\Omega_{\beta}(\overline{X},H)\big)R_1\right) +\sum_{i=|\beta_1|+1}^{k-k_0}\lambda_i\left(R_2^T\big(S(\overline{U}^T_{\beta}W\overline{V}_{\beta})-2\Omega_{\beta}(\overline{X},H)\big)R_2\right). 
\end{eqnarray*} 
Therefore, together with \eqref{eq:eta-WinT2} and \eqref{eq:Delta=case1}, we obtain that for any $(W,\eta)\in{\cal T}^2(H)$, $\Xi(W,\eta)\leq 0$. Also, it is easy to check that there exists $(W^{*},\eta^{*})\in{\cal T}^2(H)$ such $\Xi(W^{*},\eta^{*})=0$. 
%(e.g., $W^{*}\in\Re^{m\times n}$ such that ${\cal B}(W^*)=2{\cal B}(H)({\cal B}(\overline{X})-\overline{\sigma}_{k}I_{m+n})^{\dag}{\cal B}(H)$ and $\eta^*=\sum_{i=1}^{k}\sigma''_i(\overline{X};H,W^*)$). 

By combining \eqref{eq:support-product} and \eqref{eq:Xi=0-1}, we obtain that
\begin{eqnarray}
	\delta^{*}_{{\cal T}^2(H)}(\overline{S},-1) &=& \sup\left\{\langle \overline{S},W\rangle-\eta\mid (W,\eta)\in{\cal T}^2(H) \right\} \nonumber\\ [3pt]
	 &=& \sum_{l=1}^{r_0}{\rm tr}\left(2\Omega_{a_l}(\overline{X},H)\right)+\left\langle {\rm Diag}(\overline{u}_\beta), 2\Omega_{\beta}(\overline{X},H)\right\rangle. \label{eq:delta-1}
\end{eqnarray}

{\bf Case 2.} $\overline{\sigma}_k=0$. Let $\alpha$ and $\beta$ be the index sets defined by \eqref{eq:def-alpha-beta-case2} and $\beta_1$, $\beta_2$ and $\beta_3$ be the index sets defined by \eqref{eq:def-beta123-case1}. Let $r_{0}\in\{0,1,\ldots,r+1\}$ and $r_0\le\widetilde{r}_0\le r_0+1$ be the integers such that \eqref{eq:widetilde_r-case2} holds. For each $l\in\{1,\ldots,r_0\}$, since $\overline{\sigma}_i=\overline{\sigma}_{i'}$ for any $i,i'\in a_l$, we still use $\overline{\nu}_l$ to denote the common value. By \eqref{eq:eqiv-critical-1} and \eqref{eq:eqiv-critical-2}, we know that there exist orthogonal matrices  $E\in{\cal O}^{|\beta|}$ and $F\in{\cal O}^{|\beta|+n-m}$ such that \eqref{eq:simultaneous-ordered-2} holds, i.e., $[{\rm Diag}(\overline{u}_{\beta})\ \ 0]$ and $[\overline{U}_{\beta}^TH\overline{V}_{\beta}\
\ \overline{U}_{\beta}^TH\overline{V}_{2}]$ admit a simultaneous ordered SVD, which implies that $E$ and $F$ have the block diagonal structure \eqref{eq:R-block-diag-2}. Therefore, we know from the part (ii) and (iii) of Proposition \ref{prop:second-directional-diff-singlevalue} that  $(W,\eta)\in{\cal T}^2(H)$ if and only if $(W,\eta)$ satisfies the following conditions: if $\sigma_{k-k_0}\left([\overline{U}_{\beta}^TH\overline{V}_{\beta}\
\ \overline{U}_{\beta}^TH\overline{V}_{2}]\right)>0$, then
\begin{eqnarray}
&&\sum_{i=1}^{k}\sigma''_i(\overline{X};H,W) \nonumber \\ [3pt]
&=&\sum_{l=1}^{r_{0}}{\rm tr}\,(S(\overline{U}^T_{a_{l}}W\overline{V}^T_{a_{l}}))-2\sum_{l=1}^{r_0}{\rm tr}\left(\Omega_{a_l}(\overline{X},H)\right)\nonumber\\ 
&&+{\rm tr}\,\left(S\big(Q_1^{T}[\overline{U}_{\beta}^{T}(W-2H\,\overline{X}^{\dag}\,H)\overline{V}_{\beta}\ \ \overline{U}_{\beta}^{T}(W-2H\,\overline{X}^{\dag}\,H)\overline{V}_{2}]Q_1\big) \right)\nonumber\\ 
&&+\sum_{i=|\beta_1|+1}^{k-k_0}\lambda_i\left(S\big(Q_2^{T}[\overline{U}_{\beta}^{T}(W-2H\,\overline{X}^{\dag}\,H)\overline{V}_{\beta}\ \ \overline{U}_{\beta}^{T}(W-2H\,\overline{X}^{\dag}\,H)\overline{V}_{2}]Q_2\big) \right)\le \eta
\,;\label{eq:eta-WinT2-condition-1}
\end{eqnarray}
if $\sigma_{k-k_0}\left([\overline{U}_{\beta}^TH\overline{V}_{\beta}\
\ \overline{U}_{\beta}^TH\overline{V}_{2}]\right)=0$, then
\begin{eqnarray}
&&\sum_{i=1}^{k}\sigma''_i(\overline{X};H,W) \nonumber \\ [3pt]
&=&\sum_{l=1}^{r_{0}}{\rm tr}\,(S(\overline{U}^T_{a_{l}}W\overline{V}^T_{a_{l}}))-2\sum_{l=1}^{r_0}{\rm tr}\left(\Omega_{a_l}(\overline{X},H)\right)\nonumber\\ 
&&+{\rm tr}\,\left(S\big(Q_1^{T}[\overline{U}_{\beta}^{T}(W-2H\,\overline{X}^{\dag}\,H)\overline{V}_{\beta}\ \ \overline{U}_{\beta}^{T}(W-2H\,\overline{X}^{\dag}\,H)\overline{V}_{2}]Q_1\big) \right)\nonumber\\ 
&&+\sum_{i=|\beta_1|+1}^{k-k_0}\sigma_i\left({Q'_2}^{T}[\overline{U}_{\beta}^{T}(W-2H\,\overline{X}^{\dag}\,H)\overline{V}_{\beta}\ \ \overline{U}_{\beta}^{T}(W-2H\,\overline{X}^{\dag}\,H)\overline{V}_{2}]Q''_2 \right)\le \eta, \label{eq:eta-WinT2-condition-2}
\end{eqnarray}
where $\Omega_{a_l}(\overline{X},H)\in{\cal S}^{m}$, $l=1,\ldots,r_0$ are given by \eqref{eq:def-Omega-dd} with respect to $\overline{X}$, $Q_1\in{\cal O}^{|\beta_1|}$, $Q_2\in{\cal O}^{|\beta_2|}$,  $Q_3\in{\cal O}^{|b|}$ and $Q'_{3}\in{\cal O}^{|b|+n-m}$ are given by \eqref{eq:R-block-diag-2}, $Q'_2\in{\cal O}^{|\beta_2|+|\beta_3|}$ and $Q''_2\in{\cal O}^{|\beta_2|+|b|+n-m}$ are defined by
\[
Q'_2=\left[\begin{array}{cc}
Q_2 & 0 \\
0 & Q_3
\end{array}\right] \quad {\rm and} \quad Q''_2=\left[\begin{array}{cc}
Q_2 & 0 \\
0 & Q'_3
\end{array}\right].
\]
Meanwhile, since $\overline{u}_{\alpha}=e_{\alpha}$, we have for any $(W,\eta)\in{\cal T}^2(H)$, 
\begin{eqnarray}
&&-\eta+\langle \overline{S}, W \rangle = -\eta+ \left\langle \overline{U}^T\overline{S}\overline{V}, \overline{U}^TW\overline{V} \right\rangle \nonumber\\ [3pt]
&=& -\eta +\left\langle \left[\begin{array}{ccc}
{\rm Diag}(\overline{u}_{\alpha}) & 0 & 0 \\
0 & {\rm Diag}(\overline{u}_\beta) & 0
\end{array}\right], \left[\begin{array}{ccc}
S(\overline{U}_{\alpha}^TW\overline{V}_{\alpha}) & 0 & 0 \\
0 & \overline{U}_{\beta}^TW\overline{V}_{\beta} & \overline{U}_{\beta}^TW\overline{V}_{2}
\end{array}\right] \right\rangle \nonumber\\ [3pt]
&=& -\eta+\sum_{l=1}^{r_0}{\rm tr}\,(S(\overline{U}^T_{a_{l}}W\overline{V}^T_{a_{l}}))+\left\langle [{\rm Diag}(\overline{u}_\beta)\ \ 0]
,  [\overline{U}_{\beta}^TW\overline{V}_{\beta}\ \ \overline{U}_{\beta}^TW\overline{V}_{2}]  \right\rangle\nonumber\\ [3pt]
&=& \Xi(W,\eta)  +2\sum_{l=1}^{r_0}{\rm tr}\left(\Omega_{a_l}(\overline{X},H)\right)\nonumber\\ [3pt]
&&+\left\langle [{\rm Diag}(\overline{u}_\beta)\ \ 0],  2[\overline{U}^{T}_{\beta}H\,\overline{X}^{\dag}\,H\overline{V}_{\beta}\ \ \overline{U}^{T}_{\beta}H\,\overline{X}^{\dag}\,H\overline{V}_{2}] \right\rangle,\label{eq:support-product2}
\end{eqnarray}
where
\begin{eqnarray}
\Xi(W,\eta)&=&-\eta+ \sum_{l=1}^{r_0}{\rm tr}\,(S(\overline{U}^T_{a_{l}}W\overline{V}^T_{a_{l}}))-2\sum_{l=1}^{r_0}{\rm tr}\left(\Omega_{a_l}(\overline{X},H)\right)\nonumber\\ 
&&+\left\langle [{\rm Diag}(\overline{u}_\beta)\ \ 0],  [\overline{U}^{T}_{\beta}(W-2H\,\overline{X}^{\dag}\,H)\overline{V}_{\beta}\ \ \overline{U}^{T}_{\beta}(W-2H\,\overline{X}^{\dag}\,H)\overline{V}_{2}] \right\rangle.\label{eq:Delta=case2}
\end{eqnarray}

Similarly, we are able to show that
\begin{equation}\label{eq:Xi=0-2}
\max \left\{\Xi(W,\eta) \,|\, (W,\eta)\in{\cal T}^2(H)\right\}=0.
\end{equation}
In fact, if $\sigma_{k-k_0}\left([\overline{U}_{\beta}^TH\overline{V}_{\beta}\
\ \overline{U}_{\beta}^TH\overline{V}_{2}]\right)>0$, then since $0\le\overline{u}_{\beta}\le e_\beta$ and $\langle e_{\beta},\overline{u}_{\beta}\rangle\le 
k-k_0$,  we know from Lemma \ref{lem:Fan} (Fan's inequality)  that the last term of \eqref{eq:Delta=case2} satisfies
\begin{eqnarray*}
&&\left\langle [{\rm Diag}(\overline{u}_\beta)\ \ 0],  [\overline{U}^{T}_{\beta}(W-2H\,\overline{X}^{\dag}\,H)\overline{V}_{\beta}\ \ \overline{U}^{T}_{\beta}(W-2H\,\overline{X}^{\dag}\,H)\overline{V}_{2}] \right\rangle \nonumber\\ 
&=& \left\langle [{\rm Diag}(\overline{u}_\beta)\ \ 0],  E^T[\overline{U}^{T}_{\beta}(W-2H\,\overline{X}^{\dag}\,H)\overline{V}_{\beta}\ \ \overline{U}^{T}_{\beta}(W-2H\,\overline{X}^{\dag}\,H)\overline{V}_{2}]F \right\rangle \nonumber\\ 
&\leq&{\rm tr}\,\left(S\big(Q_1^{T}[\overline{U}_{\beta}^{T}(W-2H\,\overline{X}^{\dag}\,H)\overline{V}_{\beta}\ \ \overline{U}_{\beta}^{T}(W-2H\,\overline{X}^{\dag}\,H)\overline{V}_{2}]Q_1\big) \right)\nonumber\\ 
&&+\sum_{i=|\beta_1|+1}^{k-k_0}\lambda_i\left(S\big(Q_2^{T}[\overline{U}_{\beta}^{T}(W-2H\,\overline{X}^{\dag}\,H)\overline{V}_{\beta}\ \ \overline{U}_{\beta}^{T}(W-2H\,\overline{X}^{\dag}\,H)\overline{V}_{2}]Q_2\big) \right). 
\end{eqnarray*}
Thus, together with (\ref{eq:eta-WinT2-condition-1}) and (\ref{eq:Delta=case2}), we obtain that $\Xi(W,\eta)\leq 0$ for any $(W,\eta)\in{\cal T}^2(H)$. If $\sigma_{k-k_0}\left([\overline{U}_{\beta}^TH\overline{V}_{\beta}\
\ \overline{U}_{\beta}^TH\overline{V}_{2}]\right)=0$, then by Lemma \ref{lem:vonNeumann} (von Neumann's trace inequality), we know that
\begin{eqnarray*}
	&&\left\langle [{\rm Diag}(\overline{u}_\beta)\ \ 0],  [\overline{U}^{T}_{\beta}(W-2H\,\overline{X}^{\dag}\,H)\overline{V}_{\beta}\ \ \overline{U}^{T}_{\beta}(W-2H\,\overline{X}^{\dag}\,H)\overline{V}_{2}] \right\rangle \nonumber\\ 
	&=& \left\langle [{\rm Diag}(\overline{u}_\beta)\ \ 0],  E^T[\overline{U}^{T}_{\beta}(W-2H\,\overline{X}^{\dag}\,H)\overline{V}_{\beta}\ \ \overline{U}^{T}_{\beta}(W-2H\,\overline{X}^{\dag}\,H)\overline{V}_{2}]F \right\rangle \nonumber\\ 
	&\leq&{\rm tr}\,\left(S\big(Q_1^{T}[\overline{U}_{\beta}^{T}(W-2H\,\overline{X}^{\dag}\,H)\overline{V}_{\beta}\ \ \overline{U}_{\beta}^{T}(W-2H\,\overline{X}^{\dag}\,H)\overline{V}_{2}]Q_1\big) \right)\nonumber\\ 
	&&+\sum_{i=|\beta_1|+1}^{k-k_0}\sigma_i\left({Q'_2}^{T}[\overline{U}_{\beta}^{T}(W-2H\,\overline{X}^{\dag}\,H)\overline{V}_{\beta}\ \ \overline{U}_{\beta}^{T}(W-2H\,\overline{X}^{\dag}\,H)\overline{V}_{2}]Q''_2 \right).  
\end{eqnarray*}
Together with (\ref{eq:eta-WinT2-condition-2}) and (\ref{eq:Delta=case2}), we conclude that $\Xi(W,\eta)\leq 0$ for any $(W,\eta)\in{\cal T}^2(H)$. Moreover, it is easy to check that in both case there exists $(W^*,\eta^*)\in{\cal T}^2(H)$ such that $\Xi(\eta^*,W^*)=0$ (e.g., $W^*=2H\overline{X}^{\dag}H\in\Re^{m\times n}$ and $\eta^*=\sum_{i=1}^{k}\sigma''_i(\overline{X};H,W^*)$). 

By combining \eqref{eq:support-product2} and \eqref{eq:Xi=0-2}, we obtain that
\begin{eqnarray}
\delta^{*}_{{\cal T}^2(H)}(\overline{S},-1) &=& \sup\left\{\langle \overline{S},W\rangle-\eta\mid (W,\eta)\in{\cal T}^2 \right\} \nonumber\\ [3pt]
&=&2\sum_{l=1}^{r_0}{\rm tr}\left(\Omega_{a_l}(\overline{X},H)\right)+\left\langle {\rm Diag}(\overline{u}_\beta),  2\overline{U}^{T}_{\beta}H\,\overline{X}^{\dag}\,H\overline{V}_{\beta} \right\rangle. \label{eq:delta-2}
\end{eqnarray}

We summarize the above results on the support function $\delta^{*}_{{\cal T}^2(H)}$ of the second order tangent set ${\cal T}^2(H)$ in the following proposition.

\begin{proposition}\label{prop:support-function-second-order-tangent-set}
Let $(\overline{X},\overline{S})\in\Re^{m\times n}\times\Re^{m\times n}$ be a solution of the GE \eqref{eq:GE}, i.e., $\overline{S}\in \partial\, \theta(\overline{X})$. Let $X=\overline{X}+\overline{S}$ admit the SVD \eqref{eq:SVD-X}. Denote $\overline{\sigma}=\sigma(\overline{X})$ and $\overline{u}=\sigma(\overline{S})$. For any $H\in{\cal C}(X,\partial\,\theta(\overline{X}))$, let ${\cal T}^2(H)\in\Re^{m\times n}\times\Re$ be the second order tangent set defined by \eqref{eq:second-order-tangent-used}, and $\Omega_{a_l}(\overline{X},H)\in{\cal S}^{m}$, $l=1,\ldots,r_0$ and $\Omega_{\beta}(\overline{X},H)\in{\cal S}^{m}$ be the matrices given by \eqref{eq:def-Omega-dd} with respect to $\overline{X}$. Then, the support function of ${\cal T}^2(H)$ at $(\overline{S},-1)$ is given as follows.
\begin{description}
\item[(i)] If $\overline{\sigma}_{k}>0$, then
\begin{equation*}
\delta^{*}_{{\cal T}^2(H)}(\overline{S},-1)= \sum_{l=1}^{r_0}{\rm tr}\left(2\Omega_{a_l}(\overline{X},H)\right)+\big\langle{\rm Diag}(\overline{u}_\beta),2\Omega_{\beta}(\overline{X},H)\big\rangle.
\end{equation*}
\item[(ii)]  If  $\overline{\sigma}_{k}=0$, then
\begin{equation*}
\delta^{*}_{{\cal T}^2(H)}(\overline{S},-1)=\sum_{l=1}^{r_0}{\rm tr}\left(2\Omega_{a_l}(\overline{X},H)\right)+\big\langle {\rm Diag}(\overline{u}_\beta),  2\overline{U}^{T}_{\beta}H\,\overline{X}^{\dag}\,H\overline{V}_{\beta} \big\rangle.
\end{equation*}
\end{description}
\end{proposition}

\begin{remark}
	By \eqref{eq:second-order-tangent-used}, we know that for the given $\overline{S}\in \partial\, \theta(\overline{X})$ and $H\in{\cal C}(X,\partial\,\theta(\overline{X}))$, the second order tangent set ${\cal T}^2(H)$ is the epigraph of the closed convex function $\psi:=\theta''(\overline{X};H,\cdot):\Re^{m\times n}\to\Re$. Then, the support function of ${\cal T}^2(H)$ at $(\overline{S},-1)$ obtained in Proposition \ref{prop:support-function-second-order-tangent-set} equals to the conjugate function value of $\psi$ at $\overline{S}$, i.e., 
	\[
	\psi^*(\overline{S}):=\sup\{\langle W,\overline{S}\rangle-\psi(W) \mid W\in\Re^{m\times n} \}=\delta^{*}_{{\cal T}^2(H)}(\overline{S},-1).
	\]
\end{remark}
 
\begin{definition}\label{def:Upsilon-function-MCP}
For any given $\overline{X}\in\Re^{m\times n}$, define the function $\Upsilon_{ \overline{X}}: \partial\, \theta(\overline{X})\times \Re^{m\times n}\to\Re$ by  for any $\overline{S}\in\partial\,\theta(\overline{X})$ and $H\in\Re^{m\times n}$,
if $\overline{\sigma}_{k}>0$, then
\begin{eqnarray*}
\Upsilon_{\overline{X}}\left(\overline{S},H\right):=\sum_{l=1}^{r_0}{\rm tr}\left(2\Omega_{a_l}(\overline{X},H)\right)+\left\langle{\rm Diag}(\overline{u}_\beta),2\Omega_{\beta}(\overline{X},H)\right\rangle,
\end{eqnarray*} if $\overline{\sigma}_{k}=0$, then
\begin{eqnarray*}
\Upsilon_{\overline{X}}\left(\overline{S},H\right)&:=&\sum_{l=1}^{r_0}{\rm tr}\left(2\Omega_{a_l}(\overline{X},H)\right)+\left\langle {\rm Diag}(\overline{u}_\beta),  2\overline{U}^{T}_{\beta}H\,\overline{X}^{\dag}\,H\overline{V}_{\beta} \right\rangle,
\end{eqnarray*}
where $\overline{\sigma}=\sigma(\overline{X})$, $\overline{u}=\sigma(\overline{S})$, and $\Omega_{a_l}(\overline{X},H)\in{\cal S}^{m}$, $l=1,\ldots,r_0$ and $\Omega_{\beta}(\overline{X},H)\in{\cal S}^{m}$ are given by \eqref{eq:def-Omega-dd} with respect to $\overline{X}$.
%Note that all the orthogonal matrices $\overline{P}$, $\overline{U},\overline{V}$ and the index sets $a_{l}$, $l=1,\ldots,r_{0}$ and $\beta$ are defined with respect to $(\bar{t},\overline{X})$.
\end{definition}
 
Similarly, for the dual GE \eqref{eq:GE-dual}, by employing the similar arguments, we are able to derive the general results on the support function values corresponding to the second order tangent sets of the polar cone ${\cal K}^{\circ}$. In particular, we are interesting in the support function value of the following the special second order tangent set ${\cal T}^2_{\circ}(H)$ at $H\in{\cal C}(X;\partial\,\theta^*(\overline{S}))$, which is defined by
%\[
%{\cal T}^2_{\circ}(H):={\cal T}^2_{\circ}(H,0)={\cal T}^2_{{\cal K}^{\circ}}((\overline{S},-1);(H,0))=-{\rm epi}\,\vartheta''(\overline{S};H,\cdot)=\left\{(W,\eta)\in\Re^{m\times n}\times\Re\mid \vartheta''(\overline{S};H,W)\le -\eta \right\},
%\]
\[
{\cal T}^2_{\circ}(H):={\cal T}^2_{{\cal K}^{\circ}}((\overline{S},-1);(H,0))=\left\{
\begin{array}{ll}
-{\rm epi}\,\vartheta''(\overline{S};H,\cdot) & \mbox{if $\vartheta(\overline{S})=1$,} \\ [3pt]
\Re^{m\times n}\times\Re & \mbox{if $\vartheta(\overline{S})<1$,}
\end{array}
\right.
\]
where $\vartheta=\|\cdot\|_{(k)}^*$ is the dual norm of the Ky Fan $k$-norm. For simplicity, we omit the detail proof here.

\begin{proposition}\label{prop:support-function-second-order-tangent-set-dual}
	Let $(\overline{X},\overline{S})\in\Re^{m\times n}\times\Re^{m\times n}$ be a solution of the dual GE \eqref{eq:GE-dual}. Suppose that $X=\overline{X}+\overline{S}$ has the SVD \eqref{eq:SVD-X}. Denote $\overline{\sigma}=\sigma(\overline{X})$ and $\overline{u}=\sigma(\overline{S})$. For any $H\in{\cal C}(X;\partial\,\theta^*(\overline{S}))$, let $\Omega_{\alpha\cup\beta_1}(\overline{S},H)\in{\cal S}^{|\alpha|+|\beta_1|}$ and $\Omega_{a_l}(\overline{S},H)\in{\cal S}^{|a_l|}$, $l=\widetilde{r}_0+1,\ldots,\widetilde{r}_1$ be the matrices defined by \eqref{eq:def-Omega-dd} with respect to $\overline{S}$. Then, the support function of ${\cal T}^2_{\circ}(H)$ at $(\overline{X},\theta(\overline{X}))$  is given as follows.
	\begin{description}
		\item[(i)] If $\overline{\sigma}_{k}>0$, then
		\begin{eqnarray*}
			\delta^{*}_{{\cal T}^2_{\circ}(H)}(\overline{X},\theta(\overline{X}))&=&\sum_{l=1}^{r_0}\overline{\nu}_l{\rm tr}\left(2\left(\Omega_{\alpha\cup\beta_1}(\overline{S},H)\right)_{a_la_l}\right)+\overline{\sigma}_k{\rm tr}\left(2\left(\Omega_{\alpha\cup\beta_1}(\overline{S},H)\right)_{\beta_1\beta_1}\right)\\ 
			&&+ \overline{\sigma}_k\sum_{l=\widetilde{r}_0+1}^{\widetilde{r}_1}{\rm tr}\left(2\Omega_{a_l}(\overline{S},H) \right) +\overline{\sigma}_k{\rm tr}\left(2\overline{U}_{\beta_3}^TH\overline{S}^{\dag}H\overline{V}_{\beta_3}\right) \\ [3pt]
			&&+\left\langle{\rm Diag}(\overline{\sigma}_\gamma), 2\overline{U}_{\gamma}^TH\overline{S}^{\dag}H\overline{V}_{\gamma} \right\rangle.
		\end{eqnarray*}
		\item[(ii)]  If  $\overline{\sigma}_{k}=0$, then
		\begin{equation*}
			\delta^{*}_{{\cal T}^2_{\circ}(H)}(\overline{X},\theta(\overline{X}))=
			\sum_{l=1}^{r_0}\overline{\nu}_l{\rm tr}\left(2\left(\Omega_{\alpha\cup\beta_1}(\overline{S},H)\right)_{a_la_l}\right).
		\end{equation*}
	\end{description}
\end{proposition}

\begin{definition}\label{def:Upsilon-function-MCP-dual}
	For any given $\overline{S}\in \Re^{m\times n}$, define the function $\Upsilon^{\circ}_{ \overline{S}}: \partial\,\theta^*(\overline{S})\times \Re^{m\times n}\to\Re$ by  for any $\overline{X}\in\partial\,\theta^*(\overline{S})$ and $H\in\Re^{m\times n}$,
	if $\sigma_{k}(\overline{X})>0$, then
	\begin{eqnarray*}
		\Upsilon^{\circ}_{\overline{S}}\left(\overline{X},H\right)&:=&\sum_{l=1}^{r_0}\overline{\nu}_l{\rm tr}\left(2\left(\Omega_{\alpha\cup\beta_1}(\overline{S},H)\right)_{a_la_l}\right)+\overline{\sigma}_k{\rm tr}\left(2\left(\Omega_{\alpha\cup\beta_1}(\overline{S},H)\right)_{\beta_1\beta_1}\right)\\ 
		&&+ \overline{\sigma}_k\sum_{l=\widetilde{r}_0+1}^{\widetilde{r}_1}{\rm tr}\left(2\Omega_{a_l}(\overline{S},H) \right) +\overline{\sigma}_k{\rm tr}\left(2\overline{U}_{\beta_3}^TH\overline{S}^{\dag}H\overline{V}_{\beta_3}\right) \\ [3pt]
		&&+\left\langle{\rm Diag}(\overline{\sigma}_\gamma), 2\overline{U}_{\gamma}^TH\overline{S}^{\dag}H\overline{V}_{\gamma} \right\rangle,
	\end{eqnarray*} if $\sigma_{k}(\overline{X})=0$, then
	\begin{equation*}
		\Upsilon^{\circ}_{\overline{S}}\left(\overline{X},H\right):=\sum_{l=1}^{r_0}\overline{\nu}_l{\rm tr}\left(2\left(\Omega_{\alpha\cup\beta_1}(\overline{S},H)\right)_{a_la_l}\right), 
	\end{equation*}
	where $\overline{\sigma}=\sigma(\overline{X})$, $\overline{u}=\sigma(\overline{S})$, and $\Omega_{\alpha\cup\beta_1}(\overline{S},H)\in{\cal S}^{|\alpha|+|\beta_1|}$ and $\Omega_{a_l}(\overline{S},H)\in{\cal S}^{|a_l|}$, $l=\widetilde{r}_0+1,\ldots,\widetilde{r}_1$ are given by \eqref{eq:def-Omega-dd} with respect to $\overline{S}$.
	%Note that all the orthogonal matrices $\overline{P}$, $\overline{U},\overline{V}$ and the index sets $a_{l}$, $l=1,\ldots,r_{0}$ and $\beta$ are defined with respect to $(\bar{t},\overline{X})$.
\end{definition}

It seems that the functions $\Upsilon_{\overline{X}}$ and $\Upsilon^{\circ}_{ \overline{S}}$ are quite complicate from the definitions. However, one can easily compute the values by elementary calculations. Moreover, we have the following interesting proposition on the defined functions $\Upsilon_{\overline{X}}$ and $\Upsilon^{\circ}_{\overline{S}}$. 

\begin{proposition}\label{prop:Upsilon=0}
	Let $\overline{S}\in\partial\,\theta(\overline{X})$ (or equivalently $\overline{X}\in\partial\,\theta^*(\overline{S})$) be given. Then, for any $H\in\Re^{m\times n}$, $\Upsilon_{\overline{X}}(\overline{S},H)\le 0$, $\Upsilon^{\circ}_{\overline{S}}\left(\overline{X},H\right))\le 0$. Moreover, we have
	\[
	\Upsilon_{\overline{X}}(\overline{S},H)=0\quad \Longleftrightarrow\quad \Upsilon^{\circ}_{\overline{S}}\left(\overline{X},H\right)=0,
	\]
	which is equivalent to the following conditions.
\begin{itemize}
	\item[(i)] If $\sigma_k(\overline{X})>0$, then
	\begin{equation}\label{eq:condition-Upsilon-1}
	\left\{\begin{array}{l}
	\left[\begin{array}{ccc}
	\widetilde{H}_{\alpha\alpha} & \widetilde{H}_{\alpha\beta_{1}} & \widetilde{H}_{\alpha\beta_{2}} \\ 
	\widetilde{H}_{\beta_{1}\alpha} & \widetilde{H}_{\beta_{1}\beta_{1}} & \widetilde{H}_{\beta_{1}\beta_{2}} \\
	\widetilde{H}_{\beta_{2}\alpha} & \widetilde{H}_{\beta_{2}\beta_{1}} & \widetilde{H}_{\beta_{2}\beta_{2}} \end{array}\right]\in{\cal  S}^{|\alpha|+|\beta_{1}|+|\beta_{2}|},\\ [3pt]
	\widetilde{H}_{\beta_{1}\beta_{3}}=(\widetilde{H}_{\beta_{3}\beta_{1}})^{T},\quad \widetilde{H}_{\beta_{2}\beta_{3}}=(\widetilde{H}_{\beta_{3}\beta_{2}})^{T}\,\\ [3pt] 
	\widetilde{H}_{\alpha\beta_{2}}=(\widetilde{H}_{\beta_{2}\alpha})^{T}=0,\quad \widetilde{H}_{\alpha\beta_{3}}=(\widetilde{H}_{\beta_{3}\alpha})^{T}=0,\\  [3pt]
	\widetilde{H}_{\alpha\gamma}=(\widetilde{H}_{\gamma\alpha})^{T}=0,\\  [3pt]
	\widetilde{H}_{\beta_{1}\gamma}=(\widetilde{H}_{\gamma\beta_{1}})^{T}=0,\quad \widetilde{H}_{\beta_{2}\gamma}=(\widetilde{H}_{\gamma\beta_{2}})^{T}=0,\\  [3pt]
	\widetilde{H}_{\alpha c}=0,\quad \widetilde{H}_{\beta_{1} c}=0,\quad \widetilde{H}_{\beta_{2} c}=0,
	\end{array}
	\right.
	\end{equation} where $\widetilde{H}=\overline{U}^{T}H\overline{V}$, and the index sets $\alpha$, $\beta$, $\gamma$, and $\beta_{i}$, $i=1,2,3$ are defined by (\ref{eq:def-alpha-beta-gamma-case1}) and (\ref{eq:def-beta123-case1}).
	 
	\item[(ii)] If $\sigma_k(\overline{X})=0$, then
	\begin{equation}\label{eq:condition-Upsilon-2}
	\left\{\begin{array}{l}
	\widetilde{H}_{\alpha\alpha}\in{\cal  S}^{|\alpha|},\quad \widetilde{H}_{\alpha\beta_1}=(\widetilde{H}_{\beta_1\alpha})^T\\ 
	\widetilde{H}_{\alpha\beta_{2}}=(\widetilde{H}_{\beta_{2}\alpha})^{T}=0,\quad \widetilde{H}_{\alpha\beta_{3}}=(\widetilde{H}_{\beta_{3}\alpha})^{T}=0,\\ 
	\widetilde{H}_{\alpha c}=0,
	\end{array}
	\right.
	\end{equation}
	where $\widetilde{H}=\overline{U}^{T}H\overline{V}$, and the index sets $\alpha$, $\beta$, and $\beta_{i}$, $i=1,2,3$ are defined by (\ref{eq:def-alpha-beta-case2}) and (\ref{eq:def-beta123-case1}).
\end{itemize}
\end{proposition}
\noindent {\bf Proof.} Let $X=\overline{X}+\overline{S}$ admit the SVD \eqref{eq:SVD-X}. Denote $\overline{\sigma}=\sigma(\overline{X})$ and $\overline{u}=\sigma(\overline{S})$. Let $\widetilde{H}_1=\overline{U}^TH\overline{V}_1$ and $\widetilde{H}_2=\overline{U}^TH\overline{V}_2$. Consider the following two cases.

{\bf Case 1.} $\overline{\sigma}_k>0$. By \eqref{eq:def-Omega-dd} and the definition of the pseudoinverse, we obtain that 
\begin{eqnarray*}
	{\rm tr}\left(2 \Omega_{a_l}(\overline{X},H)\right) &=& \sum_{l'=1 \atop  l'\neq l}^{r+1}\frac{2}{\bar{\nu}_{l'}-\bar{\nu}_l}\|S(\widetilde{H}_1)_{a_la_{l'}}\|^2+\sum_{l'=1}^{r+1}\frac{2}{-\bar{\nu}_{l'}-\bar{\nu}_l}\|T(\widetilde{H}_1)_{a_la_{l'}}\|^2 \\ [3pt]
	&& + \frac{1}{-\bar{\nu}_l}\|(\widetilde{H}_2)_{a_l}\|^2, \quad l=1,\ldots,r_0
\end{eqnarray*}
and
\begin{eqnarray*}
	\bar{\mu}_l{\rm tr}\left(2\Omega_{a_l}(\overline{X},H)\right)&=& \sum_{l'=1}^{r_0}\frac{2\bar{\mu}_l}{\bar{\nu}_{l'}-\overline{\sigma}_k}\|S(\widetilde{H}_1)_{a_la_{l'}}\|^2+\sum_{l'=r_1+1}^{r+1}\frac{2\bar{\mu}_l}{\bar{\nu}_{l'}-\overline{\sigma}_k}\|S(\widetilde{H}_1)_{a_la_{l'}}\|^2 \\ [3pt]
	&&+\sum_{l'=1}^{r+1}\frac{2\bar{\mu}_l}{-\bar{\nu}_{l'}-\overline{\sigma}_k}\|T(\widetilde{H}_1)_{a_la_{l'}}\|^2 + \frac{\bar{\mu}_l}{-\overline{\sigma}_k}\|(\widetilde{H}_2)_{a_l}\|^2, \quad l=r_0+1,\ldots,\widetilde{r}_1.
\end{eqnarray*}
Thus, since $\overline{u}_i=0$ if $i\in\beta_3$, we have
\begin{equation*}
	\left\langle {\rm Diag}(\overline{u}_\beta), 2\Omega_{\beta}(\overline{X},H) \right\rangle= \sum_{l=r_0+1}^{\widetilde{r}_1} \bar{\mu}_l{\rm tr}\left(2\Omega_{a_l}(\overline{X},H)\right).
\end{equation*}
Therefore, we obtain the following explicit formula of $\Upsilon_{\overline{X}}(\overline{S},H)$:
\begin{eqnarray}
\Upsilon_{\overline{X}}(\overline{S},H) &=& \sum_{l=1}^{r_0}\sum_{l'=\widetilde{r}_0+1}^{\widetilde{r}_1}\frac{2(1-\bar{\mu}_{l'})}{\overline{\sigma}_k-\bar{\nu}_l}\|S(\widetilde{H}_1)_{a_la_{l'}}\|^2 \nonumber \\ [3pt]
&&+\sum_{l=1}^{r_0}\sum_{l'=\widetilde{r}_1+1}^{r+1}\frac{2}{\bar{\nu}_{l'}-\bar{\nu}_l}\|S(\widetilde{H}_1)_{a_la_{l'}}\|^2+\sum_{l=\widetilde{r}_0+1}^{\widetilde{r}_1}\sum_{l'=r_1+1}^{r+1}\frac{2\bar{\mu}_l}{\bar{\nu}_{l'}-\overline{\sigma}_k}\|S(\widetilde{H}_1)_{a_la_{l'}}\|^2 \nonumber\\ [3pt]
&&+\sum_{l=1}^{r_0}\sum_{l'=1}^{r+1}\frac{2}{-\bar{\nu}_{l'}-\bar{\nu}_l}\|T(\widetilde{H}_1)_{a_la_{l'}}\|^2+\sum_{l=r_0+1}^{\widetilde{r}_1}\sum_{l'=1}^{r+1}\frac{2\bar{\mu}_l}{-\bar{\nu}_{l'}-\overline{\sigma}_k}\|T(\widetilde{H}_1)_{a_la_{l'}}\|^2 \nonumber\\ [3pt]
&& + \sum_{l=1}^{r_0}\frac{-1}{\bar{\nu}_l}\|(\widetilde{H}_2)_{a_l}\|^2 + \sum_{l=r_0+1}^{\widetilde{r}_1}\frac{-\bar{\mu}_l}{\overline{\sigma}_k}\|(\widetilde{H}_2)_{a_l}\|^2. \label{eq:Upsilon-case1}
\end{eqnarray}
Since
\begin{equation}\label{eq:u-v-facts}
\left\{\begin{array}{ll}
\overline{\sigma}_k<\bar{\nu}_l, & l=1,\ldots,r_0, \\
\bar{\nu}_l<\overline{\sigma}_k, & l=r_1+1,\ldots,r+1, \\
\bar{\mu}_{l}<1, & l=\widetilde{r}_0+1,\ldots,\widetilde{r}_1, \\
\bar{\nu}_{l'}<\bar{\nu}_l, & l=1,\ldots,r_0,\ l'=\widetilde{r}_1+1,\ldots,r+1, \\
\bar{\nu}_l>0, & l=1,\ldots,\widetilde{r}_1,
\end{array}
\right.
\end{equation}
it is easy to see that all the coefficients of the quadric terms of \eqref{eq:Upsilon-case1} are negative, which implies that $\Upsilon_{\overline{X}}(\overline{S},H)\le 0$ and $\Upsilon_{\overline{X}}(\overline{S},H)=0$ if and only if $H\in\Re^{m\times n}$ satisfies the conditions \eqref{eq:condition-Upsilon-1}.

Meanwhile, by \eqref{eq:def-Omega-dd} and the pseudoinverse, we obtain that
\begin{eqnarray*}
	\bar{\nu}_l{\rm tr}\left(2\left(\Omega_{\alpha\cup\beta_1}(\overline{S},H)\right)_{a_la_l}\right)&=&\sum_{l'=\widetilde{r}_0+1}^{\widetilde{r}_1}\frac{2\bar{\nu}_l}{\bar{\mu}_{l'}-1}\|S(\widetilde{H}_1)_{a_la_{l'}}\|^2+\sum_{l'=\widetilde{r}_1+1}^{r+1}\frac{2\bar{\nu}_l}{-1}\|S(\widetilde{H}_1)_{a_la_{l'}}\|^2 \\ [3pt]
	&&+\sum_{l'=1}^{r+1}\frac{2\bar{\nu}_l}{-\bar{\mu}_{l'}-1}\|T(\widetilde{H}_1)_{a_la_{l'}}\|^2 + \frac{\bar{\nu}_l}{-1}\|(\widetilde{H}_2)_{a_l}\|^2, \quad l=1,\ldots,\widetilde{r}_0 
\end{eqnarray*}
and
\begin{eqnarray*}
	&& \overline{\sigma}_k{\rm tr}\left(2\Omega_{a_l}(\overline{S},H)\right) \\ [3pt]
	&=& \sum_{l'=1}^{\widetilde{r}_0}\frac{2\overline{\sigma}_k}{1-\bar{\mu}_{l}}\|S(\widetilde{H}_1)_{a_la_{l'}}\|^2+\sum_{l'=\widetilde{r}_0+1 \atop l'\neq l}^{\widetilde{r}_1}\frac{2\overline{\sigma}_k}{\bar{\mu}_{l'}-\bar{\mu}_{l}}\|S(\widetilde{H}_1)_{a_la_{l'}}\|^2+\sum_{l'=\widetilde{r}_1+1}^{r+1}\frac{2\overline{\sigma}_k}{-\bar{\mu}_l}\|S(\widetilde{H}_1)_{a_la_{l'}}\|^2 \\ [3pt]
	&&+\sum_{l'=1}^{r+1}\frac{2\overline{\sigma}_k}{-\bar{\mu}_{l'}-\bar{\mu}_l}\|T(\widetilde{H}_1)_{a_la_{l'}}\|^2 + \frac{\overline{\sigma}_k}{-\bar{\mu}_l}\|(\widetilde{H}_2)_{a_l}\|^2, \quad l=\widetilde{r}_0+1,\ldots,\widetilde{r}_1.
\end{eqnarray*}
Note that for any $A,B\in\Re^{p\times q}$, ${\rm tr}(A^TB)=\|\frac{A+B}{2}\|^2-\|\frac{A-B}{2}\|^2$.
Thus, we have for $l=\widetilde{r}_1+1,\ldots,r_1$,
\begin{eqnarray*}
	&& \overline{\sigma}_k{\rm tr}\left(2\overline{U}_{a_l}^TH\overline{S}^{\dag}H\overline{V}_{a_l} \right) \\ [3pt]
	&=& \sum_{l'=1}^{\widetilde{r}_0}\frac{2\overline{\sigma}_k}{1}\left(\|S(\widetilde{H}_1)_{a_la_{l'}}\|^2-\|T(\widetilde{H}_1)_{a_la_{l'}}\|^2\right)+\sum_{l'=\widetilde{r}_0+1}^{\widetilde{r}_1}\frac{2\overline{\sigma}_k}{\bar{\mu}_{l'}}\left(\|S(\widetilde{H}_1)_{a_la_{l'}}\|^2-\|T(\widetilde{H}_1)_{a_la_{l'}}\|^2\right), 
\end{eqnarray*}
and for $l=r_1+1,\ldots,r$,
\begin{eqnarray*}
	&& \bar{\nu}_l{\rm tr}\left(2\overline{U}_{a_l}^TH\overline{S}^{\dag}H\overline{V}_{a_l}\right) \\ [3pt]
	&=& \sum_{l'=1}^{\widetilde{r}_0}\frac{2\bar{\nu}_l}{1}\left(\|S(\widetilde{H}_1)_{a_la_{l'}}\|^2-\|T(\widetilde{H}_1)_{a_la_{l'}}\|^2\right)+\sum_{l'=\widetilde{r}_0+1}^{\widetilde{r}_1}\frac{2\bar{\nu}_l}{\bar{\mu}_{l'}}\left(\|S(\widetilde{H}_1)_{a_la_{l'}}\|^2-\|T(\widetilde{H}_1)_{a_la_{l'}}\|^2\right) .
\end{eqnarray*}
By noting that $\overline{\sigma}_i=0$ if $i\in b$,
we have
\begin{eqnarray*}
	&&\overline{\sigma}_k{\rm tr}\left(2\overline{U}_{\beta_3}^TH\overline{S}^{\dag}H\overline{V}_{\beta_3}\right) +\left\langle{\rm Diag}(\overline{\sigma}_\gamma), 2\overline{U}_{\gamma}^TH\overline{S}^{\dag}H\overline{V}_{\gamma} \right\rangle \\ [3pt]
	&=&\sum_{l=\widetilde{r}_1+1}^{r_1}\overline{\sigma}_k{\rm tr}\left(2\overline{U}_{a_l}^TH\overline{S}^{\dag}H\overline{V}_{a_l}\right)+\sum_{l=r_1+1}^{r}\bar{\nu}_l{\rm tr}\left(2\overline{U}_{a_l}^TH\overline{S}^{\dag}H\overline{V}_{a_l}\right).
\end{eqnarray*}
Therefore, we obtain the following explicit formula of $\Upsilon^{\circ}_{\overline{S}}\left(\overline{X},H\right)$:
\begin{eqnarray}
&&\Upsilon^{\circ}_{\overline{S}}\left(\overline{X},H\right) \nonumber \\ [3pt]
&=& \sum_{l=1}^{r_0}\sum_{l'=\widetilde{r}_0+1}^{\widetilde{r}_1} \left(\frac{2\bar{\nu}_l}{\bar{\mu}_{l'}-1}+\frac{2\overline{\sigma}_k}{1-\bar{\mu}_{l'}}\right)\|S(\widetilde{H}_1)_{a_la_{l'}}\|^2 \nonumber \\ [3pt]
&&+\sum_{l=1}^{r_0}\sum_{l'=\widetilde{r}_1+1}^{r_1}2(\overline{\sigma}_k-\bar{\nu}_l)\|S(\widetilde{H}_1)_{a_la_{l'}}\|^2+\sum_{l=1}^{r_0}\sum_{l'=r_1+1}^{r+1}\frac{2\bar{\nu}_l}{-1}\|S(\widetilde{H}_1)_{a_la_{l'}}\|^2\nonumber \\ [3 pt]
&&+\sum_{l=r_0+1}^{\widetilde{r}_1}\sum_{l'=r_1+1}^{r+1}\frac{2\overline{\sigma}_k}{-\bar{\mu}_l}\|S(\widetilde{H}_1)_{a_la_{l'}}\|^2+\sum_{l=1}^{\widetilde{r}_0}\sum_{l'=1}^{\widetilde{r}_1}\frac{2\bar{\nu}_l}{-\bar{\mu}_{l'}-1}\|T(\widetilde{H}_1)_{a_la_{l'}}\|^2\nonumber \\ [3pt]
&&+\sum_{l=1}^{\widetilde{r}_0}\sum_{l'=\widetilde{r}_1+1}^{r_1}2\left(\frac{-\bar{\nu}_l}{\bar{\mu}_{l'}+1}-\overline{\sigma}_k\right)\|T(\widetilde{H}_1)_{a_la_{l'}}\|^2+\sum_{l=1}^{\widetilde{r}_0}\sum_{l'=r_1+1}^{r}2\left(\frac{-\bar{\nu}_l}{\bar{\mu}_{l'}+1}-\bar{\nu}_{l'}\right)\|T(\widetilde{H}_1)_{a_la_{l'}}\|^2\nonumber \\ [3pt]
&&+\sum_{l=\widetilde{r}_0+1}^{\widetilde{r}_1}-2\bar{\nu}_l\|T(\widetilde{H}_1)_{a_la_{r+1}}\|^2  + \sum_{l=1}^{\widetilde{r}_0}-\bar{\nu}_l\|(\widetilde{H}_2)_{a_l}\|^2 + \sum_{l=\widetilde{r}_0+1}^{\widetilde{r}_1}\frac{\overline{\sigma}_k}{-\bar{\mu}_l}\|(\widetilde{H}_2)_{a_l}\|^2. \label{eq:Upsilon-o-case1}
\end{eqnarray}
Again, it follows from \eqref{eq:u-v-facts} that all the coefficients of the quadric terms of \eqref{eq:Upsilon-o-case1} are negative, which implies that $\Upsilon^{\circ}_{\overline{S}}\left(\overline{X},H\right)\le 0$ and $\Upsilon^{\circ}_{\overline{S}}\left(\overline{X},H\right)=0$ if and only if $H\in\Re^{m\times n}$ satisfies the conditions \eqref{eq:condition-Upsilon-1}.

{\bf Case 2.} $\overline{\sigma}_k=0$. By the similar arguments, we are able to show that for any $H\in\Re^{m\times n}$,
\begin{eqnarray}
\Upsilon_{\overline{X}}(\overline{S},H)&=&\sum_{l=1}^{r_0}\sum_{l'=\widetilde{r}_0+1}^{\widetilde{r}_1}2\frac{\bar{\mu}_{l'}-1}{\bar{\nu}_l}\|S(\widetilde{H}_1)_{a_la_{l'}}\|^2+\sum_{l=1}^{r_0}\sum_{l'=\widetilde{r}_1+1}^{r+1}\frac{2}{-\bar{\nu}_l}\|S(\widetilde{H}_1)_{a_la_{l'}}\|^2 \nonumber \\ [3pt] &&+\sum_{l=1}^{r_0}\sum_{l'=1}^{r_0}\frac{2}{-\bar{\nu}_{l'}-\bar{\nu}_l}\|T(\widetilde{H}_1)_{a_la_{l'}}\|^2+\sum_{l=1}^{r_0}\sum_{l'=r_0+1}^{\widetilde{r}_0}\frac{4}{-\bar{\nu}_l}\|T(\widetilde{H}_1)_{a_la_{l'}}\|^2 \nonumber \\ [3pt] &&+\sum_{l=1}^{r_0}\sum_{l'=\widetilde{r}_0+1}^{\widetilde{r}_1}2\frac{\bar{\mu}_l+1}{-\bar{\nu}_l}\|T(\widetilde{H}_1)_{a_la_{l'}}\|^2 +\sum_{l=1}^{r_0}\sum_{l'=\widetilde{r}_1+1}^{r+1}\frac{2}{-\bar{\nu}_l}\|T(\widetilde{H}_1)_{a_la_{l'}}\|^2 \nonumber \\ [3pt]
&&+\sum_{l=1}^{r_0} \frac{1}{-\bar{\nu}_l}\|(\widetilde{H}_2)_{a_l}\|^2 \label{eq:Upsilon-case2}
\end{eqnarray}
and
\begin{eqnarray}
\Upsilon^{\circ}_{\overline{S}}\left(\overline{X},H\right)&=& \sum_{l=1}^{r_0}\sum_{l'=\widetilde{r}_0+1}^{\widetilde{r}_1}\frac{2\bar{\nu}_l}{\bar{\mu}_{l'}-1}\|S(\widetilde{H}_1)_{a_la_{l'}}\|^2+\sum_{l=1}^{r_0}\sum_{l'=\widetilde{r}_1+1}^{r+1}\frac{2\bar{\nu}_l}{-1}\|S(\widetilde{H}_1)_{a_la_{l'}}\|^2 \nonumber\\ [3pt]
&&+\sum_{l=1}^{r_0}\sum_{l'=1}^{r+1}\frac{2\bar{\nu}_l}{-\bar{\mu}_{l'}-1}\|T(\widetilde{H}_1)_{a_la_{l'}}\|^2 + \sum_{l=1}^{r_0}\frac{\bar{\nu}_l}{-1}\|(\widetilde{H}_2)_{a_l}\|^2. \label{eq:Upsilon-o-case2}
\end{eqnarray} 
Thus, it follows from the fact \eqref{eq:u-v-facts} that all the coefficients of the quadric terms of \eqref{eq:Upsilon-case2} and \eqref{eq:Upsilon-o-case2} are negative, which implies that both $\Upsilon_{\overline{X}}(\overline{S},H)\le 0$, $\Upsilon^{\circ}_{\overline{S}}\left(\overline{X},H\right))\le 0$ and 
\[
\Upsilon_{\overline{X}}(\overline{S},H)=0\quad \Longleftrightarrow\quad \Upsilon^{\circ}_{\overline{S}}\left(\overline{X},H\right)=0,
\]
which is equivalent to the conditions \eqref{eq:condition-Upsilon-2}.
$\hfill\Box$

\section{Conclusions}
In this paper, we studied some important variational properties of the Ky Fan $k$-norm of matrices related to the nonlinear optimization problem involving the Ky Fan $k$-norm, which frequently arises and plays a crucial role in various applications. In particular, we introduced and study the concepts of nondegeneracy, strict complementary and critical cone to the locally optimal solutions of the basic nonlinear optimization model \eqref{eq:P-g}. Moreover, we provide the explicit formula of the conjugate function of the parabolic second order directional derivative of the Ky Fan $k$-norm, which provides the necessary second order information for the study of nonlinear optimization problem involving the Ky Fan $k$-norm. The variational results obtained in this paper can be applied immediately to the study of various perturbation and sensitivity properties, e.g., the second order optimality conditions, strong regularity, full stability and calmness of the general Ky Fan $k$-norm related optimization problems.

%\begin{acknowledgements}
%If you'd like to thank anyone, place your comments here
%and remove the percent signs.
%\end{acknowledgements}

% BibTeX users please use one of
%\bibliographystyle{spbasic}      % basic style, author-year citations
%\bibliographystyle{spmpsci}      % mathematics and physical sciences
%\bibliographystyle{spphys}       % APS-like style for physics
%\bibliography{}   % name your BibTeX data base

\begin{thebibliography}{999}

\bibitem{Bhatia97}
Bhatia, R.: Matrix Analysis. Springer Science \& Business Media (1997).

%\bibitem{BCShapiro98}
%Bonnans, J.F., Cominetti, R., Shapiro, A.: Sensitivity Analysis of Optimization Problems Under Second Order Regular Constraints. Mathematics of Operations Research 23, 806--831 (1998).

%\bibitem{BCShapiro99}
%Bonnans, J.F., Cominetti, R., Shapiro, A.: Second order optimality conditions based on parabolic second order tangent sets. SIAM Journal on Optimization 9, 466--492 (1999).
	
%\bibitem{BShapiro98}
%Bonnans, J.F., Shapiro, A.: Optimization Problems with Perturbations: A Guided Tour. SIAM Review 40, 228--264 (1998).

\bibitem{BShapiro00} 
{Bonnans, J.F., Shapiro, A.}:, Perturbation Analysis of
Optimization Problems, Springer, New York, 2000.

\bibitem{BDXiao04}
Boyd, S., Diaconis, P., Xiao, L.: Fastest Mixing Markov Chain on a Graph. SIAM Review 46, 667--689 (2004).

\bibitem{CanRec08}
Cand\`{e}s, E., Recht, B.: Exact Matrix Completion via Convex Optimization. Foundations of Computational Mathematics 9, 717--772 (2009).
 
 
\bibitem{CanTao09}
Cand\`{e}s, E., Tao, T.: The Power of Convex Relaxation: Near-Optimal Matrix Completion. IEEE Transactions on Information Theory 56, 2053--2080 (2010).

\bibitem{CLMW09}
Cand\`{e}s, E., Li, X., Ma, Y., Wright, J.:
{Robust principal component analysis?}
Journal of the ACM 58, article No. 11 (2011).


%\bibitem{CSun08}
%Chan, Z.X., Sun, D.: Constraint Nondegeneracy, Strong Regularity, and Nonsingularity in Semidefinite Programming. SIAM Journal on Optimization 19, 370--396 (2008).

\bibitem{CSPW09}
Chandrasekaran, V., Sanghavi,  S., Parrilo, P.A., Willsky, A.:
{Rank-sparsity incoherence for matrix decomposition}. SIAM Journal of Optimization  21, 572--596 (2011).

\bibitem{CFP03}
Chu, M.,  Funderlic, R., Plemmons,R.:
{Structured low rank approximation}.
Linear Algebra and its Applications 366, 157--172 (2003).

%\bibitem{Clarke76}
%Clarke, F.H. On the inverse function theorem, Pacific J. Math, 64, 97--102 (1976).

\bibitem{Clarke83}
{Clarke, F.H.}:
Optimization and Nonsmooth Analysis, John Wiley \& Sons, New York, 1983.

%\bibitem{Ding12}
%Ding, C.: An Introduction to a Class of Matrix Optimization Problems. PhD thesis. National University of Singapore, \url{http://www.math.nus.edu.sg/~matsundf/DingChao_Thesis_final.pdf} (2012).

%\bibitem{DSSToh14}
%Ding, C., Sun, D., Sun, J., Toh, K.-C.: Spectral operators of matrices, (2014). Preprint available at \url{http://arxiv.org/abs/1401.2269}.

\bibitem{DSToh10}
Ding, C., Sun, D., Toh, K.-C.: An introduction to a class of matrix cone programming. Mathematical Programming 144, 141--179 (2014).

\bibitem{Dobrynin04}
Dobrynin, V.:
{On the rank of a matrix associated with a graph}.
Discrete Mathematics 276, 169--175 (2004).

%\bibitem{Eaves71}
%Eaves, B.C.: On the basic theorem of complementarity. Mathematical Programming 1, 68--75 (1971).

\bibitem{FPang03}
Facchinei, F., Pang, J.S., eds. Finite-dimensional variational inequalities and complementarity problems. Vol. 1. Springer Science \& Business Media, 2003.

\bibitem{Fan49}
{Fan, K.}: On a theorem of Weyl concerning eigenvalues of linear transformations.  Proceedings of the
National Academy of Sciences of U.S.A. 35, 652--655  (1949).


%\bibitem{FKoranyi94}
%Faraut, J., Kor\'anyi, A.: Analysis on symmetric cones. Oxford university press, 1994.


%\bibitem{FHMangasarian84}
%Fujiwara, O., Han, S.-P., Mangasarian, O. L.: Local duality of nonlinear programs. SIAM Jounal on Control and Optimization 22 162–-169 (1984).

\bibitem{GSun10}
Y. Gao, Y., Sun, D.:
{A majorized penalty approach for calibrating rank constrained correlation matrix problems}.
Preprint available at \url{http://www.math.nus.edu.sg/~matsundf/MajorPen_May5.pdf} (2010).

\bibitem{GT94}
{Greenbaum, A., Trefethen, L.N.}:
GMRES/CR and Arnoldi/Lanczos as matrix approximation problems.
{SIAM Journal on Scientific Computing} 15, 359--368 (1994).

\bibitem{Gro09}
Gross, D.: Recovering low-rank matrices from few coefficients in any basis. IEEE Transactions on Information Theory 57, 1548--1566 (2009).

%\bibitem{HJohnson85}
% {Horn, R.A. and Johnson, C.R.}: Matrix Analysis,
%Cambridge University Press, Cambridge,  1985.

\bibitem{KLVempala97}
Kotlov, A.,  Lov\'{a}sz, L., Vempala, S.:
{The Colin de Verdi\`{e}re number and sphere representations of a graph}.
Combinatorica 17, 483--521 (1997).

\bibitem{KMO09}
Keshavan, R.H., Montanari, A., Oh, S.: Matrix completion from noisy entries. Journal of Machine Learning Research 11, 2057–2078 (2010).

%\bibitem{Kummer91}
%Kummer, B.: Lipschitzian inverse functions, directional derivatives, and applications in $C^{1,1}$ optimization. Journal of Optimization Theory and Applications 70, 561--582 (1991).

\bibitem{KKummer02}
Klatte, D., Kummer, B.: Nonsmooth Equations in Optimization: regularity, calculus, methods, and applications. Kluwer Academic Publishers, (2002).


\bibitem{Lancaster64}
{Lancaster, P.}: On eigenvalues of matrices dependent on a parameter, {Numerische Mathematik}
6, 377--387 (1964).

%\bibitem{LOverton96}
%Lewis, A.S., Overton, M.: Eigenvalue optimization. Acta Numerica 5, 149--190 (1996).

\bibitem{LSendov05}
{Lewis, A.S., Sendov, H.S.}:
Nonsmooth analysis of singular values. Part II: applications. {Set-Valued Analysis} 13, 243--264 (2005).

\bibitem{Lovasz79}
Lov\'{a}sz, L.:
{On the Shannon capacity of a graph}.
IEEE Transactions on Information Theory 25, 1--7 (1979).

%\bibitem{MSZhao05}
%Meng, F., Sun, D., Zhao, G.: Semismoothness of solutions to generalized equations and the Moreau-Yosida regularization. Mathematical Programming 104, 561--581 (2005).

\bibitem{MSPan12}
Miao, W.M., Sun, D.,  Pan, S.:
{A rank-corrected procedure for matrix completion with
	fixed basis coefficients}.
Preprint available at \url{http://arxiv.org/abs/1210.3709} (2012).

\bibitem{Mordukhovich06a}
Mordukhovich, B.S.: Variational Analysis and Generalized Differentiation I: Basic Theory. Springer Berlin Heidelberg (2006).

%\bibitem{Mordukhovich06b}
%Mordukhovich, B.S.: Variational Analysis and Generalized Differentiation II: Applications. Springer Berlin Heidelberg (2006).

\bibitem{MRSarabi13}
Mordukhovich, B.S., Rockafellar, R.T., Sarabi, M.E.: Characterizations of Full Stability in Constrained Optimization. SIAM Journal of Optimization 23, 1810--1849 (2013).

\bibitem{Moreau65}
{Moreau, J.-J.}:
{Proximit\'{e} et dualit\'{e} dans un espace hilbertien}.
{Bulletin de la Soci\'{e}t\'{e} Math\'{e}matique de France} 93, 1067--1070 (1965).

\bibitem{vonNeumann37}
{von\ Neumann, J.}:
 Some matrix-inequalities and metrization of matric-space, Tomsk University Review 1, 286--300 (1937). In:  Collected Works, Pergamon, Oxford, 1962, Volume IV, 205-218.
 
%\bibitem{Overton88}
%Overton, M.: On minimizing the maximum eigenvalue of a symmetric matrix.
%SIAM Journal on Matrix Analysis and Applications 9, 256--268 (1988).

%\bibitem{OWomersley92}
%Overton, M., Womersley, R.: On the sum of the largest eigenvalues of a
%symmetric matrix. SIAM Journal on Matrix Analysis and Applications 13, 41--45 (1992).

\bibitem{OWomersley93}
Overton, M., Womersley, R.: Optimality conditions and duality theory for minimizing sums of the largest eigenvalues of symmetric matrices. Mathematical Programming 62, 321--357  (1993).

%\bibitem{Qi09}
%Qi, H.: Local duality of nonlinear semidefinite programming.
%Mathematics of Operations Research 34 124--141 (2009).
	
\bibitem{RFParrilo10}
Recht, B., Fazel, M., Parrilo, P.A.: Guaranteed minimum-rank solutions of linear matrix equations via nuclear norm minimization. SIAM Review 52, 471--501 (2010).



%\bibitem{Recht11}
%Recht, B.: A simpler approach to matrix completion. Journal of Machine Learning Research 12, 3413--3430 (2011).



\bibitem{Robinson76}
Robinson, S.M.: First order conditions for general nonlinear optimization. SIAM Journal on Applied Mathematics 30, 597--607 (1976).


%\bibitem{Robinson80}
%Robinson, S.M.: Strongly regular generalized equations. Mathematics of Operations Research 5, 43--62 (1980).

\bibitem{Robinson84}
Robinson, S.M.: Local structure of feasible sets in nonlinear programming. II: Nondegeneracy. Mathematical Programming Study 22, 217--230 (1984).

\bibitem{Robinson87}
Robinson, S.M.: Local structure of feasible sets in nonlinear programming. III: Stability and sensitivity. Mathematical Programming Study 30, 45--66 (1987).

\bibitem{Rockafellar70}
Rockafellar, R.T.: Convex Analysis. Princeton University Press (1970).

\bibitem{RWets98}
Rockafellar, R.T., Wets, R.J.B.: Variational Analysis. Springer Berlin Heidelberg, Berlin, Heidelberg (1998).


%\bibitem{Shapiro03}
%Shapiro, A.: Sensitivity Analysis of Generalized Equations. Journal of Mathematical Sciences. 115, 2554--2565 (2003).


\bibitem{Sun06}
Sun, D.: The Strong Second-Order Sufficient Condition and Constraint Nondegeneracy in Nonlinear Semidefinite Programming and Their Implications. Mathematics of Operations Research 31, 761--776 (2006). 

%\bibitem{SSun08}
%Sun, D., Sun, J.: L\"{o}wner's Operator and Spectral Functions in Euclidean Jordan Algebras. Mathematics of Operations Research 33, 421--445 (2008).

\bibitem{Toh97}
{Toh, K.-C.}:
GMRES vs. ideal GMRES.
{SIAM Journal on Matrix Analysis and Applications} 18, 30--36 (1997).

\bibitem{TT98}
{Toh, K.-C. and Trefethen, L.N.}:
{\em The Chebyshev polynomials of a matrix}.
{SIAM Journal on Matrix Analysis and Applications} 20, 400--419 (1998).

\bibitem{Torki01}
{Torki, M.}: Second-order directional derivatives of all eigenvalues of a symmetric matrix. {Nonlinear
 Analysis} 46, 1133--1150 (2001).

%\bibitem{Warga81}
%Warga, J.: Fat homeomorphisms and unbounded derivate containers. Journal of Mathematical Analysis and Applications 81, 545--560 (1981).

\bibitem{Watson92}
{Watson, G.A.}:
{Linear best approximation using a class of polyhedral norms}.
Numerical Algorithms 2, 321--336 (1992).

\bibitem{WMGR09}
Wright, J.,  Ma, Y., Ganesh, A., Rao, S.:
{Robust principal component analysis: exact recovery of corrupted low-rank matrices via convex optimization}.
In Y. Bengio, D. Schuurmans,
J. Lafferty and C. Williams, editors, Advances in Neural Information Processing Systems 22, (2009).

\bibitem{WDSToh11}
Wu, B., Ding, C., Sun, D., Toh, K.-C.: On the Moreau--Yosida regularization of the vector $k$-norm related functions. SIAM Journal on Optimization 24, 766--794 (2014).

\bibitem{ZZXiao13}
Zhang, L., Zhang, N., Xiao, X.: On the Second-order Directional Derivatives of Singular Values of Matrices and Symmetric Matrix-valued Functions. Set-Valued and Variational Analysis 21(3), 557--586 (2013).

\bibitem{ZKurcyusz79}
Zowe, J., Kurcyusz, S.: 
Regularity and stability for the mathematical programming problem in Banach spaces, Applied Mathematics and Optimization  5, 49--62 (1979).


\end{thebibliography}

% Non-BibTeX users please use

\end{document}